\documentclass[12pt,reqno]{amsart}
\usepackage{amsmath,amssymb,amscd,latexsym,amsthm,mathrsfs}
\usepackage{color}
\usepackage[unicode]{hyperref}
\usepackage{hypbmsec,longtable}
\ifx\pdfoutput\undefined\else
\hypersetup{
colorlinks=true,
linkcolor=blue,
citecolor=blue,
urlcolor=blue,
filecolor=blue,
bookmarksnumbered=true,
pdfstartview=FitH,
pdfhighlight=/N
}
\fi
\newtheorem{theorem}{Theorem}
\newtheorem{lemma}{Lemma}
\newtheorem*{corollary}{Corollary}
\newtheorem{proposition}{Proposition}
\theoremstyle{remark}
\newtheorem*{remark}{Remark}
\newtheorem*{acknowledgements}{Acknowledgements}
\let\eps\varepsilon
\renewcommand{\d}{{\mathrm d}}
\renewcommand{\Im}{\operatorname{Im}}
\renewcommand{\Re}{\operatorname{Re}}
\newcommand\m{{\operatorname{m}}}

\newcommand\sn{\operatorname{sn}}
\newcommand\cn{\operatorname{cn}}
\newcommand\dn{\operatorname{dn}}
\newcommand\res{\operatorname{Res}}

\begin{document}

\hypersetup{pdfauthor={Mathew Rogers, Wadim Zudilin},%
pdftitle={From $L$-series of elliptic curves to Mahler measures}}

\title{From $L$-series of elliptic curves to~Mahler~measures}

\author{Mathew Rogers}
\address{Department of Mathematics and Statistics, Universit\'e de Montr\'eal,
CP 6128 succ.\ Centre-ville, Montr\'eal Qu\'ebec H3C\,3J7, Canada}
\email{mathewrogers@gmail.com}

\author{Wadim Zudilin}
\address{School of Mathematical and Physical Sciences,
The University of Newcastle, Callaghan NSW 2308, AUSTRALIA}
\email{wadim.zudilin@newcastle.edu.au}

\thanks{The first author is supported by National Science Foundation award DMS-0803107.
The second author is supported by Australian Research Council grant DP110104419.}

\dedicatory{To the blessed memory of Philippe Flajolet}

\date{December 14, 2010. {\it Last update\/}: September 8, 2011}


\begin{abstract}
We prove the conjectural relations between Mahler measures and
$L$-values of elliptic curves of conductors~20 and~24.  We also
present new hypergeometric expressions for $L$-values of CM
elliptic curves of conductors~27 and~36. Furthermore, we prove a
new functional equation for the Mahler measure of the polynomial
family $(1+X)(1+Y)(X+Y)-\alpha XY$, $\alpha\in\mathbb R$.
\end{abstract}

\maketitle

\section{Introduction}
\label{s-intro}

The Mahler measure of a two-variate Laurent polynomial $P(X,Y)$ is defined by
$$
\m(P):=\iint_{[0,1]^2}\log|P(e^{2\pi it},e^{2\pi is})|\,\d t\,\d s.
$$
In this paper we are mostly concerned with the Mahler measures of
three polynomial families,
\begin{align*}
m(\alpha)&:=\m\biggl(\alpha+X+\frac{1}{X}+Y+\frac{1}{Y}\biggr),
\\
g(\alpha)&:=\m\bigl((1+X)(1+Y)(X+Y)-\alpha XY\bigr),
\\
n(\alpha)&:=\m\bigl(X^3+Y^3+1-\alpha XY).
\end{align*}
Based on numerical experiments, Boyd observed that these functions
can be related to the values of $L$-series of elliptic
curves~\cite{Bo1}. For example, he hypothesized that
\begin{align}
m(8)=4m(2)&=\frac{24}{\pi^2}L(E_{24},2), \label{m(8)}
\\
g(4)=\frac34n(\sqrt[3]{32}) &=\frac{10}{\pi^2}L(E_{20},2),
\label{g(4)}
\end{align}
where $E_{24}$ and $E_{20}$ are elliptic curves of conductors $24$
and $20$, respectively. The primary goal of this article is to
present rigorous proofs of \eqref{m(8)} and~\eqref{g(4)}.  In the
remainder of the introduction we briefly describe our method, define
notation, review facts about Mahler measures, and present additional
theorems.

The modularity theorem shows that
the $L$-functions of elliptic curves can be equated to Mellin
transforms of weight-two modular forms. For a generic elliptic
curve $E$, we can write
\begin{equation}
L(E,2)=-\int_{0}^{1}f(q)\,\log q\,\frac{\d q}{q},
\label{LE2a}
\end{equation}
where $f(e^{2\pi i \tau})$ is a newform of weight~2 on a congruence
subgroup of $SL_2(\mathbb{Z})$.  The choice of $f(q)$ is dictated by
the elliptic curve~$E$.  For instance, if $E$ has conductor $20$,
then $f(q)=\eta^2(q^2)\eta^2(q^{10})$; if $E$ has conductor $24$,
then $f(q)=\eta(q^2)\eta(q^4)\eta(q^6)\eta(q^{12})$. For convenience
we consider the eta function with respect to $q$:
\begin{equation*}
\eta(q):=q^{1/24}\prod_{k=1}^{\infty}(1-q^{k})
=\sum_{n=-\infty}^{\infty}(-1)^nq^{(6n+1)^2/24}.
\end{equation*}
Our first step is to find modular functions $x(q)$, $y(q)$, and
$z(q)$ which depend on $f(q)$, such that
\begin{equation}
-\int_{0}^{1}f(q)\log q \,\frac{\d q}{q}
=\int_{0}^{1}x(q)\log y(q)\,\d z(q).
\label{LE2b}
\end{equation}
Next express $x$ and $y$ as algebraic functions of $z$. If we
write $x(q)=X(z(q))$, and $y(q)=Y(z(q))$, then the substitution
reduces $L(E,2)$ to a complicated(!) integral of elementary
functions:
\begin{equation*}
L(E,2)=\int_{z(0)}^{z(1)}X(z)\log Y(z)\,\d z.
\end{equation*}
Formulas for $L$-values of elliptic curves of conductor $27$,
$24$, and $20$ are given in equations \eqref{F(1,3) integral},
\eqref{G1c}, and \eqref{F(1,5) elementary integral}, respectively. The final
step is to relate the integrals to Mahler measures. We accomplish
this reduction by using properties of hypergeometric functions.

The only known approach \cite{Br}, \cite{Me} for reducing \eqref{LE2a}
to the form \eqref{LE2b} is as follows.
The quantity $L(E,2)$ is related to a Mahler measure, by first considering the
convolution $L$-function obtained from multiplying $L(E,1)$ and
$L(E,2)$.  The convolution $L$-function is then related to a
certain integral over the fundamental domain of~$E$, involving $f$
and two Eisenstein series of weight~2.  The integrals are
evaluated by Rankin's method, and the value of $L(E,1)$ cancels
from either side of the equation. Our method  is different and
works by decomposing the cusp form $f(e^{2\pi i\tau})$ into a
product of two weight~1 Eisenstein series.   We perform the
modular involution $\tau\mapsto-1/\tau$ in one of the Eisenstein
series, and then make a simple analytic change of variables in the
integral \eqref{LE2a}. As a result of these manipulations,
$L(E,2)$ reduces to an elementary integral. The details of this
computation are given in our proofs of Propositions~\ref{prop1},
\ref{prop3} and \ref{prop5} below.

Let us note that $m(\alpha)$, $n(\alpha)$ and $g(\alpha)$ can all
be expressed in terms of hypergeometric functions. These formulas
provide an efficient way to compute the Mahler measures
numerically. It was shown by Rodriguez-Villegas \cite{RV} that for
every $\alpha\in\mathbb C$,
\begin{equation}
m(\alpha)=\Re\biggl(\log\alpha
-\frac{2}{\alpha^2}\,{}_4F_3\biggl(\begin{matrix} \frac32, \, \frac32, \, 1, \, 1 \\
2, \, 2, \, 2 \end{matrix}\biggm| \frac{16}{\alpha^2} \biggr) \biggr);
\label{hyper-m1}
\end{equation}
furthermore~\cite{KO} if $\alpha\ge 0$, then
\begin{equation}
m(\alpha)=\frac\alpha4\Re{}_3F_2\biggl(\begin{matrix} \frac12, \, \frac12, \, \frac12 \\
1, \, \frac{3}{2}\end{matrix}\biggm| \frac{\alpha^2}{16} \biggr).
\label{hyper-m2}
\end{equation}
More involved hypergeometric expressions are known for $g(\alpha)$ and $n(\alpha)$ \cite[Theorem~3.1]{Rgsubmit};
in particular, the formulas
\begin{align}
g(\alpha)
&=\frac13\Re\biggl(\log\frac{(\alpha+4)(\alpha-2)^4}{\alpha^2}
-\frac{2\alpha^2}{(\alpha+4)^3}\,{}_4F_3\biggl(\begin{matrix} \frac43, \, \frac53, \, 1, \, 1 \\
2, \, 2, \, 2 \end{matrix}\biggm| \frac{27\alpha^2}{(\alpha+4)^3} \biggr)
\notag\\ &\qquad
-\frac{8\alpha}{(\alpha-2)^3}\,{}_4F_3\biggl(\begin{matrix} \frac43, \, \frac53, \, 1, \, 1 \\
2, \, 2, \, 2 \end{matrix}\biggm| \frac{27\alpha}{(\alpha-2)^3} \biggr) \biggr)
\label{hyper-g}
\end{align}
and
\begin{equation}
n(\alpha)
=\Re\biggl(\log\alpha
-\frac{2}{\alpha^3}\,{}_4F_3\biggl(\begin{matrix} \frac43, \, \frac53, \, 1, \, 1 \\
2, \, 2, \, 2 \end{matrix}\biggm| \frac{27}{\alpha^3} \biggr) \biggr)
\label{hyper-n}
\end{equation}
are valid for $|\alpha|$ sufficiently large.  Formula
\eqref{hyper-g} can also be shown to hold on the real line if
$\alpha\in\mathbb R\setminus[-4,2]$.

It is a subtle but important point that our proofs are essentially
elementary.  The modularity theorem shows that $L(E,2)=L(f,2)$,
however the formulas we prove for $L(f,2)$ are true
\emph{unconditionally}. For example, many of Boyd's conjectures
can be restated as relations between Mahler measures and the
quadruple lattice sum~\cite{Rgsubmit}
\begin{equation}
\begin{split}
F(b,c)&:=(b+1)^2(c+1)^2
\\ &\qquad
\times\sum_{\substack{n_i=-\infty\\i=1,2,3,4}}^{\infty}
\frac{(-1)^{n_1+n_2+n_3+n_4}}{\bigl((6n_1+1)^2+b(6n_2+1)^2+c(6n_3+1)^2+bc(6n_4+1)^2\bigr)^2},
\end{split}
\label{F(b,c)}
\end{equation}
where the default method of summation, is ``summation by cubes"
\cite{Bor}.  Since $L(E_{20},2)=F(1,5)$, and $L(E_{24},2)=F(2,3)$,
the above examples can be written as
\begin{align}
m(8)=4m(2)&=\frac{24}{\pi^2}F(2,3),
\label{m(8)F}
\\
g(4)=\frac34n(\sqrt[3]{32})
&=\frac{10}{\pi^2}F(1,5).
\label{g(4)F}
\end{align}
Formulas \eqref{m(8)F} and \eqref{g(4)F} are true even without the
modularity theorem.  In fact, one significant aspect of Boyd's
work, is that it provides a recipe to relate slowly-converging
lattice sums to hypergeometric functions. For more details, many
other conjectural examples as well as for state-of-art in the
area, the reader may consult \cite{Bo1}, \cite{RV}
and~\cite{Rgsubmit}.

We prove many additional theorems with the strategy we have
described.  For instance, we construct new hypergeometric
evaluations
\begin{align*}
L(E_{27},2) &=\frac{\Gamma^3(\tfrac13)}{27}\,
\,{}_3F_2\biggl(\begin{matrix} \frac13, \, \frac13, \, 1 \\
\frac23, \, \frac43 \end{matrix}\biggm| 1 \biggr)
-\frac{\Gamma^3(\tfrac23)}{18}
\,{}_3F_2\biggl(\begin{matrix} \frac23, \, \frac23, \, 1 \\
\frac43, \, \frac53 \end{matrix}\biggm| 1 \biggr),
\\
L(E_{36},2)
&=-\frac{2\pi^2\log2}{27}+\frac{\Gamma^3(\tfrac13)}{3\cdot2^{7/3}}
\,{}_3F_2\biggl(\begin{matrix} \frac13, \, \frac13, \, 1 \\
\frac56, \, \frac43 \end{matrix}\biggm| -\frac18 \biggr)
+\frac{\Gamma^3(\tfrac23)}{2^{11/3}}
\,{}_3F_2\biggl(\begin{matrix} \frac23, \, \frac23, \, 1 \\
\frac76, \, \frac53 \end{matrix}\biggm| -\frac18 \biggr)
\end{align*}
for the $L$-series of CM elliptic curves of conductors~27 and~36. We
also prove elementary integrals for lattice sums which are not
associated to elliptic curves: $F(3,7)$, $F(6,7)$ and $F(3/2,7)$.
Finally, we derive a new functional equation
\begin{equation*}
g\bigl(4p(1+p)\bigr)+g\biggl(\frac{4(1+p)}{p^2}\biggr)=2g\biggl(\frac{2(1+p)^2}p\biggr),
\qquad \frac{\sqrt3-1}2\le p\le 1,
\end{equation*}
for the Mahler measure $g(\alpha)$.  This last formula resembles
some of the functional equations due to Lal\'\i n and Rogers~\cite{LR}.

We conclude the introduction with a word about notation. This paper
involves a large number of $q$-series manipulations, and draws
heavily from Berndt's versions of Ramanujan's Notebooks
\cite{Be3,Be4,Be5}, and from Ramanujan's Lost Notebook~\cite{BA}.
For this reason, we have chosen to preserve Ramanujan's theta
function notation
\begin{equation}\label{phipsidef}
\varphi(q):=\sum_{n=-\infty}^{\infty}q^{n^2}, \qquad
\psi(q):=\sum_{n=0}^{\infty}q^{n(n+1)/2}.
\end{equation}
We also define the notation for signature~$3$ theta functions in
the next section.

\section{Conductor 27}
\label{s-cond27}

In this section we look at the CM elliptic curves of conductors $27$
and $36$, as well as at some non-elliptic curve lattice sums. Recall
that elliptic curves of conductor $27$ are associated to
$\eta^2(q^3)\eta^2(q^9)$, and elliptic curves of conductor $36$ are
associated to $\eta^4(q^6)$ \cite{Ono}. It follows that
$L(E_{27},2)=F(1,3)$ and $L(E_{36},2)=F(1,1)$. Define
\begin{equation*}
\begin{split}
H(x):=&\int_{0}^{1}\frac{\eta^3(q^3)}{\eta(q)}\,\frac{\eta^3(q^x)}{\eta(q^{3x})}\,\log
q\,\frac{\d q}{q}\\
=&\frac{1}{3}\int_{0}^{1}b(q^{x})\,c(q)\log q\,\frac{\d q}{q},
\end{split}
\end{equation*}
where the signature-$3$ theta functions are given by
\begin{align*}
a(q):=&\sum_{m,n=-\infty}^{\infty}q^{m^2+m n+n^2},\\
b(q):=&\frac{1}{2}\bigl(3a(q^3)-a(q)\bigr)=\frac{\eta^3(q)}{\eta(q^{3})},\\
c(q):=&\frac{1}{2}\bigl(a(q^{1/3})-a(q)\bigr)=3\frac{\eta^3(q^3)}{\eta(q)}.
\end{align*}
The functions $a(q)$, $b(q)$ and $c(q)$ were studied in great detail
by Ramanujan and the Borweins \cite{Be5}, \cite{Borwein}. They form
the basis of the theory of signature~$3$ theta functions.
The following lemma shows that certain values of~\eqref{F(b,c)} are expressed in
terms of~$H(x)$.

\begin{lemma}\label{Lemma on special values of H}
The following relations are true:
\begin{align}
9L(E_{27},2)&=9F(1,3)=-H(1),
\label{F(1,3) in terms of H}\\
36L(E_{36},2)&=36 F(1,1)
=-4H\biggl(\frac{4}{3}\biggr)+\frac{1}{4}H\biggl(\frac{1}{12}\biggr),
\label{F(1,1) in terms of H}\\
\frac{27}{16}F(3,7)
&=\frac{8}{7}H(1)-H(7)-\frac{1}{49}H\biggl(\frac{1}{7}\biggr),
\label{F(3,7) in terms of H}\\
\frac{27}{49}F(6,7)
&=\frac{1}{49}H\biggl(\frac{2}{7}\biggr)+H(14)-\frac{8}{7}H(2),
\label{F(6,7) in terms of H}\\
\frac{27}{25}F\biggl(\frac{3}{2},7\biggr)
&=\frac{2}{7}H\biggl(\frac{1}{2}\biggr)-\frac{1}{4}H\biggl(\frac{7}{2}\biggr)-\frac{1}{14^2}H\biggl(\frac{1}{14}\biggr).
\label{F(3/2,7) in terms of H}
\end{align}
\end{lemma}

\begin{proof}
Equation \eqref{F(1,3) in terms of H} follows from the definition
of $H(x)$.  Formula \eqref{F(1,1) in terms of H} follows from
integrating a modular equation equivalent to Somos \cite[Entry
$t_{36,9,39}$]{so2}:
\begin{equation*}
3\eta^4(q^6)=-b(q)c(q^{12})+b(q^4)c(q^3).
\end{equation*}
We can recover \eqref{F(3,7) in terms of H} by integrating a
modular equation equivalent to Ramanujan \cite[pg.~236, Entry
68]{Be4}:
\begin{equation*}
9\eta(q)\eta(q^3)\eta(q^7)\eta(q^{21})=-b(q)c(q)-7b(q^7)c(q^7)+b(q^7)c(q)+b(q)c(q^7).
\end{equation*}
Equation \eqref{F(3/2,7) in terms of H} follows from a modular
equation equivalent to Somos \cite[Entry $x_{42,8,56}$]{so2}:
\begin{equation*}
9 \eta(q^2) \eta(q^3) \eta(q^{14})
\eta(q^{21})=b(q)c(q^{14})+b(q^7)c(q^2)-b(q)c(q^2)-7b(q^7)c(q^{14}),
\end{equation*}
and \eqref{F(6,7) in terms of H} follows from a modular equation
equivalent to Somos \cite[Entry $x_{42,8,64}$]{so2}:
\begin{equation*}
9 \eta(q)\eta(q^6)\eta(q^7)\eta(q^{42})=-b(q^2)c(q^7)-b(q^{14})c(q)
+b(q^2)c(q)+7b(q^{14})c(q^7).
\end{equation*}
\vskip-\baselineskip
\end{proof}

Next we prove a second integral for $H(x)$ which
involves signature~$3$ theta functions.  This is the fundamental result needed to relate values of $H(x)$ to elementary integrals.

\begin{proposition}
\label{prop1}
Suppose that $x>0$, then
\begin{equation}\label{H(x) reduced}
\begin{split}
H(x)=\frac{2\pi}{\sqrt{3}x}\int_{0}^{1}b(q)c(q^3)\log\biggl(3\frac{c(q^{9x})}{c(q^{3x})}\biggr)\frac{\d q}{q}.
\end{split}
\end{equation}
\end{proposition}

\begin{proof}
Begin by setting $q=e^{-2\pi u}$, then
\begin{equation*}
H(x)=-\frac{(2\pi)^2}{3}\int_{0}^{\infty}u\,b(e^{-2\pi x
u})\,c(e^{-2\pi u})\,\d u.
\end{equation*}
Since $b(q)=\eta^3(q)/\eta(q^3)$ and $c(q)=3\eta^3(q^3)/\eta(q)$, it follows that $b(q)$ and $c(q)$ are linked by an involution:
\begin{equation*}
c\left(e^{-2\pi/(3u)}\right)=\sqrt{3} u\, b\left(e^{-2\pi u}\right).
\end{equation*}
We will use the following Eisenstein series expansion \cite[pg.~406]{BA}:
\begin{equation*}
c(q)=3\sum_{n=1}^{\infty}\chi_{-3}(n)\biggl(\frac{q^{n/3}}{1-q^{n/3}}-\frac{q^n}{1-q^n}\biggr).
\end{equation*}

Rearranging the series, and then applying the involution, we find that
\begin{align}
c\left(e^{-2\pi u}\right)
&=3\sum_{n,k=1}^{\infty}\chi_{-3}(n)(e^{-2\pi n k u/3}-e^{-2\pi n k u}),
\\
b\left(e^{-2\pi x u}\right)
&=\frac{\sqrt{3}}{xu}\sum_{r,s=1}^{\infty}\chi_{-3}(r)(e^{-2\pi r s/(9 x u)}-e^{-2\pi r s/(3 x u)}).
\end{align}
Therefore, the integral becomes
\begin{equation*}
\begin{split}
H(x)=-\frac{(2\pi)^2\sqrt{3}}{x}\sum_{n,k,r,s\ge 1}\chi_{-3}(nr)
&\int_{0}^{\infty}(e^{-2\pi n k u/3}-e^{-2\pi n ku})
\\ &\quad\times
(e^{-2\pi r s/(9x u)}-e^{-2\pi r s/(3x u)})\, \d u.
\end{split}
\end{equation*}
Use linearity and a $u$-substitution, to regroup the integral:
\begin{equation*}
\begin{split}
H(x)=-\frac{(2\pi)^2\sqrt{3}}{x}\sum_{n,k,r,s\ge 1}\chi_{-3}(nr)
&\int_{0}^{\infty}e^{-2\pi n k u}(e^{-2\pi r s/(3xu)}
\\ &\quad
-4e^{-2\pi r s/(9x u)}+3e^{-2\pi r s/(27x u)})\, \d u.
\end{split}
\end{equation*}
Finally make the $u$-substitution $u\mapsto r u/k$. This
permutes the indices of summation inside the integral and we obtain
\begin{equation*}
\begin{split}
H(x)=-\frac{(2\pi)^2\sqrt{3}}{x}\sum_{n,k,r,s\ge 1}\frac{r\chi_{-3}(rn)}{k}
&\int_{0}^{\infty}e^{-2\pi n r u} \bigl(e^{-2\pi k s/(3xu)}
\\ &\quad
-4e^{-2\pi k s/(9x u)}+3e^{-2\pi k s/(27x u)}\bigr)\, \d u.
\end{split}
\end{equation*}
Simplifying reduces things to
\begin{equation*}
\begin{split}
H(x)=-\frac{(2\pi)^2\sqrt{3}}{x}
&\int_{0}^{\infty}\biggl(\sum_{n,r=1}^{\infty}r\chi_{-3}(r n)e^{-2\pi r n u}\biggr)
\\ &\quad\times
\log\prod_{s=1}^{\infty}\frac{(1-e^{-2\pi s/(9 x u)})^4}
{(1-e^{-2\pi s/(27x u)})^3(1-e^{-2\pi s/(3 x u)})}\,\d u.
\end{split}
\end{equation*}
Notice that the product equals a ratio of Dedekind eta functions
where all of the $q^{1/24}$ terms have cancelled out. Applying the
involution for the eta function, we obtain
\begin{equation*}
\begin{split}
H(x)=-\frac{(2\pi)^2\sqrt{3}}{x}
&\int_{0}^{\infty}\biggl(\sum_{n,r=1}^{\infty}r\chi_{-3}(r n)e^{-2\pi r n u}\biggr)
\\ &\quad\times
\log\biggr(\frac{e^{4\pi xu}}{3}
\prod_{s=1}^{\infty}\frac{(1-e^{-2\pi s(9 x u)})^4}{(1-e^{-2\pi s(27x u)})^3(1-e^{-2\pi s(3 x u)})}\biggr) \d u.
\end{split}
\end{equation*}
Set $q=e^{-2\pi u}$, and then use the product expansion
$c(q)=3\eta^3(q^3)/\eta(q)$ \cite[pg.~109]{Be5}, to obtain
\begin{equation}\label{H(x) almost done}
H(x)=\frac{(2\pi)\sqrt{3}}{x}\int_{0}^{1}\biggl(\sum_{n,r=1}^{\infty}r\chi_{-3}(r n)q^{r n}\biggr)
\log\biggl(3\frac{c(q^{9x})}{c(q^{3x})}\biggr) \frac{\d q}{q}.
\end{equation}
To simplify the Eisenstein series, notice that
\begin{equation*}
\chi_{-3}(n)=\frac{2}{\sqrt{3}}\Im(e^{2\pi i n/3}),
\end{equation*}
and therefore
\begin{equation*}
\sum_{n,r=1}^{\infty}r\chi_{-3}(r n)q^{r n}
=-\frac{1}{12\sqrt{3}}\Im L(e^{2\pi i/3}q),
\end{equation*}
where
\begin{equation}
L(q):=1-24\sum_{n=1}^{\infty}\frac{n q^n}{1-q^n}.
\label{L(q)}
\end{equation}
By Ramanujan's Eisenstein series for $a^2(q)$ \cite[pg.~100]{Be5},
we have
\begin{equation*}
2a^2(q)=3L(q^3)-L(q),
\end{equation*}
so it follows that
\begin{equation*}
\sum_{n,r=1}^{\infty}r\chi_{-3}(r n)q^{r n}
=\frac{1}{6\sqrt{3}}\Im a^2(e^{2\pi i/3} q).
\end{equation*}
Finally, if we use
\begin{equation*}
a(e^{2\pi i/3} q)=b(q)+i\sqrt{3}c(q^3),
\end{equation*}
then
\begin{equation*}
\Im a^2(e^{2\pi i/3} q)=2\sqrt{3}b(q)c(q^3),
\end{equation*}
which implies
\begin{equation}\label{H(x) eisenstein series real reduction}
\sum_{n,r=1}^{\infty}r\chi_{-3}(r n)q^{r n}=\frac{1}{3}b(q)c(q^3).
\end{equation}
Substituting \eqref{H(x) eisenstein series real reduction} into
\eqref{H(x) almost done} concludes the proof of~\eqref{H(x) reduced}.
\end{proof}

In the next proposition, we pass from an integral involving modular
functions, to a purely elementary integral.  In order to accomplish
this, we use the inversion formulas for signature~$3$ theta
functions.

\begin{proposition}\label{Theorem on H(x) in elementary integrals}
Suppose that $x>0$, and assume that $\beta$ has degree $3x$ over
$\alpha$ in the theory of signature~$3$. Then
\begin{equation}\label{H reduction integral 1}
H(x)=\frac{2\pi}{3\sqrt{3}x}
\int_{0}^{1}\frac{(1-\alpha)^{1/3}\bigl(1-(1-\alpha)^{1/3}\bigr)}{\alpha(1-\alpha)}
\,\log\frac{1-(1-\beta)^{1/3}}{\beta^{1/3}}\,\d\alpha.
\end{equation}

Now suppose that $\beta$ has degree $x$ over $\alpha$ in the
theory of signature~$3$.  Then
\begin{equation}\label{H reduction integral 2}
H(x)=\frac{2\pi}{3\sqrt{3}x}
\int_{0}^{1}\frac{\alpha^{1/3}(1-\alpha^{1/3})}{\alpha(1-\alpha)}
\,\log\frac{1-(1-\beta)^{1/3}}{\beta^{1/3}}\,\d\alpha.
\end{equation}
\end{proposition}

\begin{proof}
Let us prove \eqref{H reduction integral 1} first.  By
formulas (2.8) and (2.9) in~\cite[pg.~93--94]{Be5}, we know that
\begin{equation*}
3c(q^3)=a(q)-b(q).
\end{equation*}
Therefore \eqref{H(x) reduced} reduces to
\begin{equation*}
H(x)=\frac{2\pi}{3\sqrt{3}x}\int_{0}^{1}b(q)\bigl(a(q)-b(q)\bigr)
\log\frac{a(q^{3x})-b(q^{3x})}{c(q^{3x})}\,\frac{\d q}{q}.
\end{equation*}
Now set
$$
q=\exp\biggl(\frac{-2\pi}{\sqrt{3}}\frac{{}_2F_1(\frac{1}{3},\frac{2}{3};1;1-\alpha)}
{{}_2F_1(\frac{1}{3},\frac{2}{3};1;\alpha)}\biggr)
$$
and notice that
$$
a^2(q)\frac{\d q}{q}=\frac{\d\alpha}{\alpha(1-\alpha)}.
$$
It is also known~\cite[pg.~103]{Be5} that
$b(q)/a(q)=(1-\alpha)^{1/3}$ and $c(q)/a(q)=\alpha^{1/3}$.
Substituting these relations completes the
proof of \eqref{H reduction integral 1}.  Equation \eqref{H
reduction integral 2} follows if we first let $q\mapsto q^{1/3}$
in \eqref{H(x) reduced}, then use
\begin{equation*}
b(q^{1/3})=a(q)-c(q),
\end{equation*}
and finally make the same substitution for~$q$.
\end{proof}

While it is known that algebraic relations exist between $\alpha$
and $\beta$ for all rational values of $x$, it is very difficult
to apply those relations except in a few cases.

\begin{theorem}
We have
\begin{equation}\label{F(1,3) in terms of 3F2's}
L(E_{27},2)=\frac{\Gamma^3(\tfrac13)}{27}\,
\,{}_3F_2\biggl(\begin{matrix} \frac13, \, \frac13, \, 1 \\
\frac23, \, \frac43 \end{matrix}\biggm| 1 \biggr)
-\frac{\Gamma^3(\tfrac23)}{18}
\,{}_3F_2\biggl(\begin{matrix} \frac23, \, \frac23, \, 1 \\
\frac43, \, \frac53 \end{matrix}\biggm| 1 \biggr).
\end{equation}
\end{theorem}

\begin{proof}
If $x=1$ in \eqref{H reduction integral 2}, then $\alpha=\beta$,
and we obtain a formula for $L(E_{27},2)$:
\begin{equation}\label{F(1,3) integral}
L(E_{27},2)=-\frac{2\pi}{27\sqrt{3}}\int_{0}^{1}\frac{\alpha^{1/3}(1-\alpha^{1/3})}{\alpha(1-\alpha)}
\,\log\frac{1-(1-\alpha)^{1/3}}{\alpha^{1/3}}\,\d\alpha.
\end{equation}
It is possible to simplify \eqref{F(1,3) integral} with
\texttt{Mathematica}.  The easiest method is to make the
substitution
\begin{equation*}
\log\frac{1-(1-\alpha)^{1/3}}{\alpha^{1/3}}=\sum_{n=1}^{\infty}\frac{(1-\alpha)^{n}-3(1-\alpha)^{n/3}}{3n},
\end{equation*}
and then perform term-by-term integration using beta integrals.
\end{proof}

The new formula for $F(1,3)$ should be
compared to the well-known $_4F_3$ evaluation~\cite[Eq.~(43)]{Rgsubmit}:
\begin{equation}\label{F(1,3) 4F3 formula}
\frac{81}{4\pi^2}L(E_{27},2)=\log6 +\frac{1}{108}
\,{}_4F_3\biggl(\begin{matrix} \frac43, \, \frac53, \, 1, \, 1 \\
2, \, 2, \, 2 \end{matrix}\biggm| -\frac{1}{8} \biggr).
\end{equation}
It seems to be a tricky task to demonstrate the equivalence of
\eqref{F(1,3) in terms of 3F2's} and \eqref{F(1,3) 4F3 formula} by
purely hypergeometric techniques.

Note that an identity can be derived for $H(1/3)$, by setting $x=1/3$
in~\eqref{H reduction integral 1}.

\section{Conductor 24}
\label{s-cond24}

It is known that an elliptic curve $E_{24}$ of conductor $24$ is
associated to the eta product $\eta(q^2)\eta(q^4) \eta(q^6) \eta(q^{12})$ \cite{Ono}.
Thus $L(E_{24},2)=F(2,3)$, where $F(b,c)$ is the four-dimensional
lattice sum~\eqref{F(b,c)}.  Let us define $G(x)$ as follows:
\begin{align*}
G(x):=&\int_{0}^{1}\frac{\eta^2(q^2)}{\eta(q)}\,\frac{\eta^2(q^6)}{\eta(q^3)}
\,\frac{\eta^2(q^x)}{\eta(q^{2x})}\,\frac{\eta^2(q^{3x})}{\eta(q^{6x})}\,\log q\,\frac{\d q}{q},\\
=&\int_{0}^{1}q^{1/2}\,\psi(q)\,\psi(q^3)\,\varphi(-q^x)\,\varphi(-q^{3x})\,\log q\,\frac{\d q}{q}.
\end{align*}
The second identity is a consequence of the product expansions
\begin{align*}
q^{1/8}\psi(q)=\frac{\eta^2(q^2)}{\eta(q)},&&\varphi(-q)=\frac{\eta^2(q)}{\eta(q^2)},
\end{align*}
where $\psi(q)$ and $\varphi(q)$ are defined in \eqref{phipsidef}.
It is easy to see that $G(1)=-4L(E_{24},2)$.  It follows that we
can solve Boyd's conductor $24$ conjectures by reducing $G(1)$ to
hypergeometric functions.

\begin{proposition}
\label{prop3}
Let $\omega=e^{2\pi i/3}$.  The following formulas hold for $x>0$:
\begin{align}\label{G(x) imaginary integral}
G(x)
&=\frac{2\pi}{3x}\Im\int_{0}^{1}\omega q\psi^4(\omega^2 q^2)
\log\biggl(4q^{3x}\frac{\psi^4(q^{12x})}{\psi^4(q^{6x})}\biggr)\frac{\d q}{q}
\\
\label{G(x) real integral}
&=\frac{\pi}{2\sqrt{3}x}\int_{0}^{1}(A-B)(A-3B)(A^2-3B^2)
\log\biggl(4q^{3x/2}\frac{\psi^4(q^{6x})}{\psi^4(q^{3x})}\biggr)\frac{\d q}{q},
\end{align}
where $A=q^{1/8}\psi(q)$ and $B=q^{9/8}\psi(q^{9})$.
\end{proposition}

\begin{proof}
Begin by setting $q=e^{-2\pi u}$; then the integral becomes
\begin{equation*}
G(x)=-(2\pi)^2\int_{0}^{\infty}u\,e^{-\pi u}\,\psi(e^{-2\pi
u})\psi(e^{-6\pi u})\,\varphi(-e^{-2\pi x u})\varphi(-e^{6\pi x
u})\,\d u.
\end{equation*}
Now consider a Lambert series due to Ramanujan \cite[pg.~223, Entry~3.1]{Be3}:
\begin{equation*}
q^{1/2}\psi(q)\psi(q^3)=\sum_{n=1}^{\infty}\frac{\chi(n)q^{n/2}}{1-q^n},
\end{equation*}
where $\chi(n)$ has conductor $6$, with $\chi(5)=-1$. Rearranging
Ramanujan's result, and then using the involution for the eta
function, we have
\begin{align}
e^{-\pi u}\psi(e^{-2\pi u})\,\psi(e^{-6\pi u})&=\sum_{n,k=1}^{\infty}\chi(n)(e^{-\pi n k u}-e^{-2\pi n k u}),
\\
\varphi(-e^{-2\pi x u})\,\varphi(-e^{-6\pi x u})&=\frac{2}{\sqrt{3} x u}\sum_{r,s=1}^{\infty}\chi(r)(e^{-2\pi r s/(12x u)}-e^{-2\pi r s/(6x u)}).
\end{align}
Noting that $\chi(n)$ is totally multiplicative, the integral
becomes
\begin{align*}
G(x)
=-\frac{8\pi^2}{\sqrt{3} x}\sum_{n,k,r,s\ge 1}\chi(r n)
&\int_{0}^{\infty}(e^{-\pi n k u}-e^{-2\pi n k u})
(e^{-2\pi r s/(12x u)}-e^{-2\pi r s/(6x u)})\,\d u
\\
=-\frac{8\pi^2}{\sqrt{3} x}\sum_{n,k,r,s\ge 1}\chi(r n)
&\int_{0}^{\infty}e^{-2\pi n k u}(2e^{-2\pi r s/(24x u)}
\\ &\qquad
-3e^{-2\pi r s/(12x u)}+e^{-2\pi r s/(6x u)})\,\d u.
\end{align*}
Now make the substitution $u\mapsto r u/k$.  This step is
crucially important, because it groups the $r$ and $n$ indices
together:
\begin{align*}
G(x)=-\frac{8\pi^2}{\sqrt{3} x}\sum_{n,k,r,s\ge 1}\frac{r\chi(r n)}{k}
&\int_{0}^{\infty}e^{-2\pi r n u}(2e^{-2\pi k s/(24x u)}
\\ &\qquad
-3e^{-2\pi k s/(12x u)}+e^{-2\pi k s/(6x u)})\,\d u.
\end{align*}
Simplifying the $k$ and $s$ sums, brings the integral to
\begin{align*}
G(x)
=-\frac{8\pi^2}{\sqrt{3} x}\int_{0}^{\infty}
&\biggl(\sum_{n,r\ge1}r\chi(r n)e^{-2\pi r nu}\biggr)
\\ &\times
\log\prod_{s=1}^{\infty}\frac{(1-e^{-2\pi s/(12 x u)})^3}{(1-e^{-2\pi s/(24 x u)})^2(1-e^{-2\pi s/(6 x u)})}\,\d u.
\end{align*}
The product equals a ratio of eta functions (the $q^{1/24}$ terms
have cancelled out). Applying the involution again, we have
\begin{align*}
G(x)=-\frac{8\pi^2}{\sqrt{3} x}\int_{0}^{\infty}
&\biggl(\sum_{n,r\ge1}r\chi(r n)e^{-2\pi r nu}\biggr)
\\ &\times
\log\biggl(\frac{e^{-3\pi u x/2}}{\sqrt{2}}
\prod_{s=1}^{\infty}\frac{(1-e^{-24\pi s x u})^3}{(1-e^{-48\pi s x u})^2(1-e^{-12\pi s x u})}\biggr)\d u.
\end{align*}
Now use the product expansion $q^{1/8}\psi(q)=\eta^2(q^2)/\eta(q)$, and simplify:
\begin{align}
G(x)
&=-\frac{4\pi}{\sqrt{3} x}\int_{0}^{1}\biggl(\sum_{n,r\ge 1}r\chi(r n)q^{r n}\biggr)
\log\biggl(\frac{q^{-3 x/4}}{\sqrt{2}}\frac{\psi(q^{6x})}{\psi(q^{12x})}\biggr)
\frac{\d q}{q}
\notag\\
&=\frac{\pi}{\sqrt{3} x}\int_{0}^{1}\biggl(\sum_{n,r\ge 1}r\chi(r n)q^{r n}\biggr)
\log\biggl(4 q^{3x}\frac{\psi^4(q^{12x})}{\psi^4(q^{6x})}\biggr)\frac{\d q}{q}.
\label{G(x) almost done}
\end{align}
The calculation is nearly complete.  To simplify the Eisenstein
series, we use
\begin{equation*}
\chi(n)=\frac{1}{\sqrt{3}}\Im(e^{2\pi i n/3}-(-1)^n e^{2\pi i n/3}),
\end{equation*}
and therefore
\begin{equation*}
\sum_{n,r\ge 1}r\chi(r n)q^{r n}
=-\frac{1}{24\sqrt{3}}\Im\bigl(L(e^{2\pi i/3}q)-L(-e^{2\pi i/3} q)\bigr),
\end{equation*}
where $L(q)$ is the Eisenstein series~\eqref{L(q)}.
Ramanujan proved \cite[pg.~114, Entry~8.2]{Be3} that
\begin{equation*}
3\varphi^4(q)=4L(q^4)-L(q),
\end{equation*}
hence
\begin{equation*}
\sum_{n,r\ge 1}r\chi(r n)q^{r n}
=\frac{1}{8\sqrt{3}}\Im\bigl(\varphi^4(e^{2\pi i/3}q)-\varphi^4(-e^{2\pi i/3} q)\bigr);
\end{equation*}
finally by \cite[pg.~40]{Be3}, we have
\begin{equation}\label{eisenstein reduction}
\sum_{n,r\ge 1}r\chi(r n)q^{r n}
=\frac{2}{\sqrt{3}}\Im\bigl(e^{2\pi i/3}q\psi^4(e^{4\pi i/3}q^2)\bigr).
\end{equation}
Substituting \eqref{eisenstein reduction} into \eqref{G(x) almost
done} completes the proof of \eqref{G(x) imaginary integral}.  To
reduce \eqref{G(x) imaginary integral} to \eqref{G(x) real
integral}, we can substitute the following identity into~\eqref{G(x) imaginary integral}:
\begin{equation*}
2\psi(\omega^2q^2)=2\psi(q^2)-3q^2\psi(q^{18})-i\sqrt{3}q^2\psi(q^{18}).
\end{equation*}
\vskip-\baselineskip
\end{proof}

\begin{lemma}
\label{cor-G1}
We have
\begin{equation}
-4L(E_{24},2)=G(1)
=\frac{\pi}{12}\int_0^{1/2}\frac{\sqrt{(1-2p)(2-p)}\log\dfrac{p^3(2-p)}{1-2p}}{(1-p^2)\sqrt
p}\,\d p. \label{G1c}
\end{equation}
\end{lemma}

\begin{proof}
Set $x=1$ and then manipulate \eqref{G(x) imaginary integral}, to obtain
\begin{equation*}
G(1)=\frac{\pi}{6}\Im\int_{0}^{1} \omega q^{1/2}\psi^4(\omega^2q)
\log\biggl(16q^{3}\frac{\psi^8(q^{6})}{\psi^8(q^{3})}\biggr)\frac{\d q}{q}.
\end{equation*}
Now apply complex conjugation, then use $\omega^2=-e^{\pi i/3}=-\omega^{1/2}$,
and let $\omega q\mapsto q$, to arrive at
\begin{equation*}
G(1)= \frac{\pi}{6}\Im \int_{0}^{\omega} q^{1/2}\psi^4(q)
\log\biggl(16q^{3}\frac{\psi^8(q^{6})}{\psi^8(q^{3})}\biggr)\frac{\d q}{q}.
\end{equation*}
Now set $\alpha(q):=1-\varphi^4(-q)/\varphi^4(q)$, and
$z(q):=\varphi^2(q)$. Then by formula \cite[pg.~123, Entry
11.1]{Be3} and \cite[pg.~120, Entry 9.1]{Be3},
\begin{align*}
q^{1/2}\psi^4(q)&=\frac{\sqrt{\alpha(q)}}{4}z^2(q),
\\
\frac{\d \alpha(q)}{\d q}&=\frac{\alpha(q)(1-\alpha(q))z^2(q)}{q}.
\end{align*}
By formulas \cite[pg.~123, Entry 11.1]{Be3} and \cite[pg~123,
Entry 11.3]{Be3} we also have
\begin{equation*}
16q^{3}\frac{\psi^8(q^6)}{\psi^8(q^3)}=\alpha(q^3).
\end{equation*}
Thus,
\begin{align*}
G(1)&=\frac{\pi}{24}\Im\int_0^{\omega}\frac{\log\alpha(q^{3})}{\sqrt{\alpha(q)}(1-\alpha(q))}\,\d\alpha(q)
\\
&=\frac{\pi}{24}\Im\int_0^{1}\frac{\log\alpha(q^{3})}{\sqrt{\alpha(\omega q)}(1-\alpha(\omega q))}\,\d\alpha(\omega q).
\end{align*}
Note that both $\alpha(\omega q)$ and $\alpha(q^{3})$ vary from 0 to~1 as $q$~changes
in the range from~0 to~1, and that the path for the latter is purely real.

The functions $\alpha(q)$ and $\beta(q)=\alpha(q^3)$ are related by the modular polynomial
$$
(\alpha^2+\beta^2+6\alpha\beta)^2-16\alpha\beta\bigl(4(1+\alpha\beta)-3(\alpha+\beta)\bigr)^2=0
$$
and admit the rational parametrization
\begin{equation}
\alpha=\frac{p(2+p)^3}{(1+2p)^3}, \quad \beta=\frac{p^3(2+p)}{1+2p}
\label{param}
\end{equation}
with $p$ ranging from 0 to 1 as as $q$~changes in the range.
The same modular relation and parametrization, of course, remain true when we take
$\omega q$ for~$q$, except that in this case the parameter $p$ ranges along the
complex curve
$$
\mathcal P=\biggl\{p:0<\frac{p^3(2+p)}{1+2p}<1\biggr\}
$$
in the upper half-plane $\Im p>0$ joining the points 0 and $-1$.
This gives rise to writing $G(1)$ as
\begin{equation*}
G(1)=\frac{\pi}{24}\Im\int_{\mathcal P}\frac{\log\beta}{\sqrt\alpha(1-\alpha)}\,\d\alpha.
\end{equation*}

First note that the integrand, as function of $p$, is analytic in the half-plane
$\Im p>0$, so that we can change the path of integration to the straight interval
from~0 to~$-1$ understood as the interval along the upper cut of the real axis:
\begin{equation*}
G(1)=\frac{\pi}{24}\Im\int_0^{-1}\frac{\log\beta}{\sqrt\alpha(1-\alpha)}\,\d\alpha
=\frac{\pi}{24}\int_0^{-1}\Im\biggl(\frac{\log\beta}{\sqrt\alpha(1-\alpha)}\biggr)\d\alpha.
\end{equation*}
Secondly, along the interval $-1<p<-1/2$ the integrand is purely \emph{real},
so that
\begin{equation*}
G(1)=\frac{\pi}{24}\int_0^{-1/2}\Im\biggl(\frac{\log\beta}{\sqrt\alpha(1-\alpha)}\biggr)\d\alpha.
\end{equation*}
Developing now the substitution \eqref{param}, computing the imaginary part
and putting $-p$ for~$p$, we thus arrive at~\eqref{G1c}.
\end{proof}

\begin{remark}
A similar recipe expresses $G(1/2)$ in the form
\begin{equation}
G(1/2)=\frac\pi{24}\Im\int_{\mathcal P}\frac{\log\beta}{1-\alpha}\,\d\alpha
\label{G12}
\end{equation}
for the path~$\mathcal P$ given above. The substitution \eqref{param}
produces an expression whose anti-derivative could be expressed in
terms of the logarithmic and dilogarithmic functions, and we finally
arrive at
$$
G(1/2)=-\frac{\pi^2\log2}3.
$$
\end{remark}

\subsection{The hypergeometric reduction}

In \eqref{G1c}, $G(1)$ splits into two integrals of the form
$$
F_1(\lambda)
=\int_0^{1/\lambda}\frac{\sqrt{(1-\lambda p)(\lambda-p)}\log(1/p)}{(1-p^2)\sqrt p}\,\d p
$$
and
$$
F_2(\lambda)
=\int_0^{1/\lambda}\frac{\sqrt{(1-\lambda p)(\lambda-p)}
\log\dfrac{\lambda-p}{1-\lambda p}}{(1-p^2)\sqrt p}\,\d p,
$$
where $\lambda=2$.

\begin{lemma}
\label{lemk1}
The identity
\begin{equation*}
F_1(\lambda)-F_2(\lambda)
=\pi\cdot{}_3F_2\biggl(\begin{matrix} \frac12, \, \frac12, \, \frac12 \\ \frac32, \, 1 \end{matrix}\biggm|
\frac1{\lambda^2}\biggr)
\end{equation*}
is true for all $\lambda\ge1$.
\end{lemma}

\begin{proof}
Making the change $\hat p=(1-\lambda p)/(\lambda-p)$ in the integral defining $F_2(\lambda)$
we obtain $p=(1-\lambda\hat p)/(\lambda-\hat p)$ and
$$
F_2(\lambda)
=(\lambda^2-1)\int_0^{1/\lambda}\frac{\sqrt{\hat p}\log(1/\hat p)}
{(1-\hat p^2)\sqrt{(1-\lambda\hat p)(\lambda-\hat p)}}\,\d\hat p,
$$
Then we set $z=1/\lambda^2$ and perform the changes $p=t\sqrt z$ and $\hat p=t\sqrt z$,
so that the required identity becomes equivalent to
\begin{align}
&
(1-z)\int_0^1\frac{\sqrt t\log(t\sqrt z)}{(1-zt^2)\sqrt{(1-t)(1-zt)}}\,\d t
-\int_0^1\frac{\sqrt{(1-t)(1-zt)}\log(t\sqrt z)}{(1-zt^2)\sqrt t}\,\d t
\notag\\ &\quad
=\pi\cdot{}_3F_2\biggl(\begin{matrix} \frac12, \, \frac12, \, \frac12 \\ \frac32, \, 1 \end{matrix}\biggm| z\biggr)
\label{identity-a}
\end{align}
for $0\le z\le1$. The left-hand side here is
\begin{align}
&
\int_0^1\frac{\bigl((1-z)t-(1-t)(1-zt)\bigr)\log(t\sqrt z)}{(1-zt^2)\sqrt{t(1-t)(1-zt)}}\,\d t
=\!\int_0^1\frac{\bigl((1-zt^2)-2(1-t)\bigr)\log(t\sqrt z)}{(1-zt^2)\sqrt{t(1-t)(1-zt)}}\,\d t
\notag\\ &\quad
=\int_0^1\frac{\log(t\sqrt z)}{\sqrt{t(1-t)(1-zt)}}\,\d t
-2\int_0^1\frac{\sqrt{1-t}\log(t\sqrt z)}{(1-zt^2)\sqrt{t(1-zt)}}\,\d t.
\label{lhs-id}
\end{align}
Our strategy is to write the series expansions of
\begin{align*}
G_\eps(z)
&=\int_0^1\frac{t^\eps\,\d t}{\sqrt{t(1-t)(1-zt)}}
=\frac{\Gamma(\frac12)\Gamma(\frac12+\eps)}{\Gamma(1+\eps)}
\,{}_2F_1\biggl(\begin{matrix} \frac12, \, \frac12+\eps \\ 1+\eps \end{matrix}\biggm| z\biggr)
\\
&=\sum_{n=0}^\infty\frac{\Gamma(n+\frac12)\Gamma(n+\frac12+\eps)}{\Gamma(n+1)\Gamma(n+1+\eps)}z^n
\end{align*}
and
$$
\tilde G_\eps(z)=\sum_{n=0}^\infty g_nz^n=\int_0^1\frac{t^\eps\sqrt{1-t}}{(1-zt^2)\sqrt{t(1-zt)}}\,\d t.
$$
Because
$$
\frac1{\sqrt{1-zt}}=\sum_{k=0}^\infty\frac{(\frac12)_k}{k!}t^kz^k,
\qquad
\frac1{1-zt^2}=\sum_{m=0}^\infty t^{2m}z^m,
$$
we have
\begin{align*}
g_n
&=\sum_{k=0}^n\frac{(\frac12)_k}{k!}\int_0^1t^{2n-k-1/2+\eps}(1-t)^{1/2}\,\d t
=\sum_{k=0}^n\frac{(\frac12)_k}{k!}\,\frac{\Gamma(2n-k+\frac12+\eps)\Gamma(\frac32)}{\Gamma(2n-k+2+\eps)}
\\
&=\frac{\Gamma(\frac32)\Gamma(2n+\frac12+\eps)}{\Gamma(2n+2+\eps)}
\sum_{k=0}^n\frac{(\frac12)_k}{k!}\,\frac{(-2n-1-\eps)_k}{(-2n+\frac12-\eps)_k}
\\
\intertext{(we apply \cite[(2.6.3)]{Slater} to the partial sum of the ${}_2F_1$ series to $n+1$ terms)}
&=\frac{\Gamma(\frac32)\Gamma(2n+\frac12+\eps)}{\Gamma(2n+2+\eps)}
\,\frac{\Gamma(n+\frac32)\Gamma(-n-\eps)}{\Gamma(n+1)\Gamma(-n+\frac12-\eps)}
\\ &\qquad\times
{}_3F_2\biggl(\begin{matrix} \frac12, \, -2n-1-\eps, \, -n+\frac12-\eps \\ -n+\frac12-\eps, \, -2n+\frac12-\eps \end{matrix}\biggm| 1\biggr)
\displaybreak[2]\\
&=\frac{\Gamma(\frac32)\Gamma(2n+\frac12+\eps)}{\Gamma(2n+2+\eps)}
\,\frac{\Gamma(n+\frac32)\Gamma(-n-\eps)}{\Gamma(n+1)\Gamma(-n+\frac12-\eps)}
\,{}_2F_1\biggl(\begin{matrix} \frac12, \, -2n-1-\eps \\ -2n+\frac12-\eps \end{matrix}\biggm| 1\biggr)
\\
\intertext{(we apply the Gauss summation to the ${}_2F_1$ series)}
&=\frac{\Gamma(\frac32)\Gamma(2n+\frac12+\eps)}{\Gamma(2n+2+\eps)}
\,\frac{\Gamma(n+\frac32)\Gamma(-n-\eps)}{\Gamma(n+1)\Gamma(-n+\frac12-\eps)}
\,\frac{\Gamma(-2n+\frac12-\eps)\Gamma(1)}{\Gamma(-2n-\eps)\Gamma(\frac32)}
\\
\intertext{(finally we use the functional equations for the Gamma function)}
&=\frac{\Gamma(n+\frac32)\Gamma(n+\frac12+\eps)}{(2n+1+\eps)\Gamma(n+1)\Gamma(n+1+\eps)}.
\end{align*}
Therefore,
\begin{align*}
G_\eps(z)-2\tilde G_\eps(z)
&=\sum_{n=0}^\infty\biggl(\frac{\Gamma(n+\frac12)\Gamma(n+\frac12+\eps)}{\Gamma(n+1)\Gamma(n+1+\eps)}
\\ &\qquad
-\frac{\Gamma(n+\frac32)\Gamma(n+\frac12+\eps)}{(n+\frac12+\frac12\eps)\Gamma(n+1)\Gamma(n+1+\eps)}\biggr)z^n
\\
&=\sum_{n=0}^\infty\frac{\Gamma(n+\frac12)\Gamma(n+\frac12+\eps)}{\Gamma(n+1)\Gamma(n+1+\eps)}
\biggl(1-\frac{n+\frac12}{n+\frac12+\frac12\eps}\biggr)z^n
\\
&=\eps\sum_{n=0}^\infty\frac{\Gamma(n+\frac12)\Gamma(n+\frac12+\eps)}{\Gamma(n+1)\Gamma(n+1+\eps)(2n+1+\eps)}z^n.
\end{align*}
This implies for~\eqref{lhs-id} that
$$
\int_0^1\frac{\bigl((1-zt^2)-2(1-t)\bigr)\log(\sqrt z)}{(1-zt^2)\sqrt{t(1-t)(1-zt)}}\,\d t
=\log\sqrt z\cdot\bigl(G_\eps(z)-2\tilde G_\eps(z)\bigr)\big|_{\eps=0}
=0
$$
and
\begin{align*}
\int_0^1\frac{\bigl((1-zt^2)-2(1-t)\bigr)\log t}{(1-zt^2)\sqrt{t(1-t)(1-zt)}}\,\d t
&=\frac{\d}{\d\eps}\bigl(G_\eps(z)-2\tilde G_\eps(z)\bigr)\bigg|_{\eps=0}
\\
&=\sum_{n=0}^\infty\frac{\Gamma(n+\frac12)^2}{\Gamma(n+1)^2(2n+1)}z^n
\\
&=\Gamma(\tfrac12)^2\sum_{n=0}^\infty\frac{(\frac12)_n^2}{n!^2(2n+1)}z^n
\\
&=\pi\cdot{}_3F_2\biggl(\begin{matrix} \frac12, \, \frac12, \, \frac12 \\ \frac32, \, 1 \end{matrix}\biggm|
z\biggr),
\end{align*}
thus establishing the required identity~\eqref{identity-a}.
\end{proof}

The method also allows us to give closed forms individually for $F_1(\lambda)$ and $F_2(\lambda)$.

\begin{lemma}
\label{lemk2}
For $\lambda\ge1$,
\begin{equation}
F_1(\lambda)
=\frac\pi2\,\log(4\lambda)
+\frac\pi2\cdot{}_3F_2\biggl(\begin{matrix} \frac12, \, \frac12, \, \frac12 \\
\frac32, \, 1 \end{matrix}\biggm| \frac1{\lambda^2}\biggr)
-\frac{\pi}{16\lambda^2}\cdot{}_4F_3\biggl(\begin{matrix} \frac32, \, \frac32, \, 1, \, 1 \\
2, \, 2, \, 2 \end{matrix}\biggm| \frac1{\lambda^2}\biggr).
\label{F1}
\end{equation}
\end{lemma}

\begin{proof}
As we have shown in the proof of Lemma~\ref{lemk1}
\begin{align*}
F_1(1/\sqrt z)
&=-\int_0^1\frac{(1-t)(1-zt)\log(t\sqrt z)}{(1-zt^2)\sqrt{t(1-t)(1-zt)}}\,\d t
\\
&=-\int_0^1\frac{\bigl((1-t)-zt(1-t)\bigr)\log(t\sqrt z)}{(1-zt^2)\sqrt{t(1-t)(1-zt)}}\,\d t,
\end{align*}
and this integral can be computed by examining the constant and linear terms in
the $\eps$-expansion of
\begin{align*}
\tilde G_\eps(z)-z\tilde G_{1+\eps}(z)
&=\sum_{n=0}^\infty\frac{\Gamma(n+\frac32)\Gamma(n+\frac12+\eps)}{(2n+1+\eps)\Gamma(n+1)\Gamma(n+1+\eps)}z^n
\\ &\quad
-\sum_{n=0}^\infty\frac{\Gamma(n+\frac32)\Gamma(n+\frac32+\eps)}{(2n+2+\eps)\Gamma(n+1)\Gamma(n+2+\eps)}z^{n+1}
\\ &\kern-20mm
=\frac{\Gamma(\frac32)\Gamma(\frac12+\eps)}{\Gamma(2+\eps)}
+\sum_{n=1}^\infty\frac{\Gamma(n+\frac12)\Gamma(n+\frac12+\eps)}{\Gamma(n+1)\Gamma(n+1+\eps)}
\biggl(\frac{n+\frac12}{2n+1+\eps}-\frac{n}{2n+\eps}\biggr)z^n
\\ &\kern-20mm
=\frac{\Gamma(\frac32)\Gamma(\frac12+\eps)}{\Gamma(2+\eps)}
+\frac\eps2\sum_{n=1}^\infty\frac{\Gamma(n+\frac12)\Gamma(n+\frac12+\eps)}
{\Gamma(n+1)\Gamma(n+1+\eps)(2n+\eps)(2n+1+\eps)}z^n.
\end{align*}
Then
\begin{align*}
\int_0^1\frac{(1-t)(1-zt)\log(\sqrt z)}{(1-zt^2)\sqrt{t(1-t)(1-zt)}}\,\d t
&=\log\sqrt z\cdot\bigl(\tilde G_\eps(z)-z\tilde G_{1+\eps}(z)\bigr)\big|_{\eps=0}
\\
&=\frac{\Gamma(\frac32)\Gamma(\frac12)}{\Gamma(2)}\log\sqrt z
=\frac{\pi\log\sqrt z}2
\end{align*}
and
\begin{align*}
&
\int_0^1\frac{(1-t)(1-zt)\log t}{(1-zt^2)\sqrt{t(1-t)(1-zt)}}\,\d t
=\frac{\d}{\d\eps}\bigl(\tilde G_\eps(z)-z\tilde G_{1+\eps}(z)\bigr)\bigg|_{\eps=0}
\\ &\quad
=\frac\pi2\biggl(\frac{\Gamma'(\frac12)}{\Gamma(\frac12)}-1+\gamma\biggr)
+\frac12\sum_{n=1}^\infty\frac{\Gamma(n+\frac12)^2}
{\Gamma(n+1)^2(2n)(2n+1)}z^n
\\ &\quad
=\frac\pi2(-2\log2-1)
+\frac\pi2\sum_{n=1}^\infty\frac{(\frac12)_n^2}
{n!^2}\biggl(\frac1{2n}-\frac1{2n+1}\biggr)z^n
\\ &\quad
=-\pi\log2
+\frac{\pi z}{16}\cdot{}_4F_3\biggl(\begin{matrix} \frac32, \, \frac32, \, 1, \, 1 \\
2, \, 2, \, 2 \end{matrix}\biggm| z\biggr)
-\frac\pi2\cdot{}_3F_2\biggl(\begin{matrix} \frac12, \, \frac12, \, \frac12 \\
\frac32, \, 1 \end{matrix}\biggm| z\biggr).
\end{align*}
Joining the latter results we obtain \eqref{F1}.
\end{proof}

Using Lemmas~\ref{cor-G1}, \ref{lemk1}, \ref{lemk2},
the equality $m(8)=4m(2)$ as well as the hypergeometric
evaluations~\eqref{hyper-m1} and \eqref{hyper-m2} of $m(8)$ and $m(2)$ we finally arrive at

\begin{theorem}
\label{th-F(2,3)}
The following evaluation is true:
$$
L(E_{24},2)=-\frac14G(1)=\frac{\pi^2}6m(2).
$$
\end{theorem}

\subsection{The elliptic reduction}
In this subsection we give an alternative derivation of
Theorem~\ref{th-F(2,3)}. In order to accomplish this, we use
properties of the Jacobian elliptic functions. Recall that $\sn u$
depends implicitly on $\alpha$, and that it is doubly periodic, with
periods $4K$ and $2iK'$, where
\begin{equation*}
K=\frac{\pi}{2}\,{}_2F_1\biggl(\begin{matrix} \frac12, \, \frac12 \\
1 \end{matrix}\biggm| \alpha \biggr),
\qquad
K'=\frac{\pi}{2}\,{}_2F_1\biggl(\begin{matrix} \frac12, \, \frac12 \\
1 \end{matrix}\biggm| 1-\alpha \biggr).
\end{equation*}
We also take the usual definition of the elliptic nome, namely
\begin{equation*}
q=\exp\biggl(-\pi\frac{K'}{K}\biggr)
=\exp\biggl(-\pi\frac{{}_2F_1(\frac{1}{2},\frac{1}{2};1;1-\alpha)}
{{}_2F_1(\frac{1}{2},\frac{1}{2};1;\alpha)}\biggr)
\end{equation*}
In the first lemma we give a Fourier series expansion for a ratio
of Jacobian elliptic functions.  Formula
\eqref{ellipticexpansioncombined} is absent from most references,
however it can be derived from results in~\cite{wanli}.

\begin{lemma}
The following identity is true:
\begin{equation}\label{ellipticexpansioncombined}
\begin{split}
\frac{\cn^2u\,\dn^2u}{1-\alpha\sn^4u}
&=\frac{\pi}{4K}+\frac{\pi}{K}\sum_{n=1}^{\infty}\frac{q^n}{1+q^{2n}}\,\cos\frac{2\pi n u}{K}
\\ &\qquad
+\frac{\pi}{\sqrt{\alpha}K}\sum_{n=0}^{\infty}\frac{q^{n+1/2}}{1+q^{2n+1}}\,\cos\frac{\pi (2n+1)u}{K}.
\end{split}
\end{equation}
\end{lemma}

\begin{proof}
Equation \eqref{elliptic expansion one}
is a superposition of results in \cite{wanli}.  Let us begin by
decomposing the function using partial fractions
\begin{align*}
\frac{\cn^2u\,\dn^2u}{1-\alpha\sn^4u}
&=\frac{(1-\sn^2u)(1-\alpha\sn^2u)}{1-\alpha\sn^4u}
\\
&=-1-\frac{(1-\sqrt{\alpha})^2}{2\sqrt{\alpha}}\,\frac{1}{1-\sqrt{\alpha}\sn^2u}
\\ &\qquad
+\frac{(1+\sqrt{\alpha})^2}{2\sqrt{\alpha}}\,\frac{1}{1+\sqrt{\alpha}\sn^2u}.
\end{align*}
By equation (1.1) in \cite[~pg. 543]{wanli}, we can show that
\begin{align}
\frac{1}{1-\sqrt{\alpha}\sn^2u}
&=\frac{\Pi(\sqrt{\alpha},\alpha)}{K}
+\frac{\pi}{(1-\sqrt{\alpha})K}\sum_{n=1}^{\infty}\frac{(-1)^nq^{n/2}}{1+q^n}\,\cos\frac{\pi n u}{K},
\label{elliptic expansion one}
\\
\frac{1}{1+\sqrt{\alpha}\sn^2u}
&=\frac{\Pi(-\sqrt{\alpha},\alpha)}{K}
+\frac{\pi}{(1+\sqrt{\alpha})K}\sum_{n=1}^{\infty}\frac{q^{n/2}}{1+q^n}\,\cos\frac{\pi n u}{K},
\label{elliptic expansion two}
\end{align}
where $\Pi(\alpha,\beta)$ is the complete elliptic integral of the
third kind.  Substituting \eqref{elliptic expansion one} and
\eqref{elliptic expansion two}, we obtain
\begin{equation*}
\begin{split}
\frac{\cn^2u\,\dn^2u}{1-\alpha\sn^4u}
&=\frac{h(\alpha)}{K}+\frac{\pi}{K}\sum_{n=1}^{\infty}\frac{q^n}{1+q^{2n}}\,\cos\frac{2\pi n u}{K}
\\ &\qquad
+\frac{\pi}{\sqrt{\alpha}K}\sum_{n=0}^{\infty}\frac{q^{n+1/2}}{1+q^{2n+1}}\,\cos\frac{\pi(2n+1)u}{K},
\end{split}
\end{equation*}
where
\begin{equation*}
h(\alpha):=-K-\frac{(1-\sqrt{\alpha})^2}{2\sqrt{\alpha}}\,\Pi(\sqrt{\alpha},\alpha)
+\frac{(1+\sqrt{\alpha})^2}{2\sqrt{\alpha}}\,\Pi(-\sqrt{\alpha},\alpha).
\end{equation*}
Finally, we are grateful to James Wan for pointing out that a more
general formula for $\Pi(m,n)$ implies that $h(\alpha)=\pi/4$
(see~\cite{We}). We leave this final calculation as an exercise for
the reader.
\end{proof}

\begin{proposition}\label{conductor 24 mahler
measure integral proposition}
Suppose that $0\le \alpha \le 1$.
The following identities are true:
\begin{align}
-\frac{8} {\pi}\int_{0}^{1}\frac{\sqrt{(1-v^2)(1-\alpha v^2)}}{1-\alpha v^4}\log v\,\d v
&=m\biggl(\frac{4}{\sqrt{\alpha}}\biggr)+\frac{1}{\sqrt{\alpha}}m(4\sqrt{\alpha})
\notag\\ &\qquad
+\log\sqrt{\alpha},
\label{sn mahler measure reduction} \displaybreak[2]\\
-\frac{8}{\pi}\int_{0}^{1}\frac{\sqrt{(1-v^2)(1-\alpha v^2)}}{1-\alpha v^4}\log(1-v^2)\,\d v
&=2m\biggl(\frac{4}{\sqrt{\alpha}}\biggr)+\log\frac{\alpha}{1-\alpha}
\notag\\ &\qquad
+\frac{1}{\sqrt{\alpha}}\,\log\frac{1-\sqrt{\alpha}}{1+\sqrt{\alpha}},
\label{cn mahler measure reduction} \displaybreak[2]\\
-\frac{8}{\pi}\int_{0}^{1}\frac{\sqrt{(1-v^2)(1-\alpha v^2)}}{1-\alpha v^4}\log(1-\alpha v^2)\,\d v
&=\frac{2}{\sqrt{\alpha}}m(4\sqrt{\alpha})-\log(1-\alpha)
\notag\\ &\qquad
+\frac{1}{\sqrt{\alpha}}\,\log\frac{1-\sqrt{\alpha}}{1+\sqrt{\alpha}}.
\label{dn mahler measure redution}
\end{align}
\end{proposition}

\begin{proof}
First notice that if we set $v=\sn u$, then
\eqref{sn mahler measure reduction} becomes
\begin{equation*}
\int_{0}^{1}\frac{\sqrt{(1-v^2)(1-\alpha v^2)}}{1-\alpha v^4}\log v\,\d v
=\int_{0}^{K}\frac{\cn^2u\,\dn^2u}{1-\alpha\sn^4u}\,\log\sn u\,\d u.
\end{equation*}
We now substitute Fourier expansions for Jacobian elliptic
functions.  The following series holds for $u\in(0,K)$
\cite[pg.~917]{GR}:
\begin{equation}\label{log sn expansion}
\log\sn u=\log\frac{2K}\pi+\log\sin\frac{\pi u}{2K}
-2\sum_{n=1}^{\infty}\frac{1}{n}\frac{q^n}{1+q^n}\biggl(1-\cos\frac{\pi n u}{K}\biggr).
\end{equation}
Substitute \eqref{ellipticexpansioncombined} and \eqref{log sn
expansion} into the integral, and then integrate term-by-term.  It
is necessary to use the following formula several times:
\begin{equation*}
\int_{0}^{K}\cos\frac{\pi n u}{K}\,\log\sin\frac{\pi u}{2K}\,\d u=\begin{cases}
-K\log 2 &\text{if $n=0$},\\
-K/(2n) &\text{if $n\ge 1$}.
\end{cases}
\end{equation*}
A substantial amount of work reduces the integral to
\begin{equation*}
\begin{split}
\int_{0}^{K}\frac{\cn^2u\,\dn^2u}{1-\alpha\sn^4u}\,\log\sn u\,\d u
&=\frac{\pi}{4}\biggl(\log\frac K\pi-2\sum_{n=1}^{\infty}\frac{1}{n}\,\frac{q^n}{1+q^n}\biggr)
\\ &\qquad
+\frac{\pi}{2}\sum_{n=1}^{\infty}\frac{1}{n}\,\frac{q^n}{1+q^{2n}}\biggl(\frac{q^{2n}}{1+q^{2n}}-\frac{1}{2}\biggr)
\\ &\qquad
+\frac{\pi}{\sqrt{\alpha}}\sum_{n=0}^{\infty}\frac{1}{2n+1}\,\frac{q^{n+1/2}}{1+q^{2n+1}}
\biggl(\frac{q^{2n+1}}{1+q^{2n+1}}-\frac{1}{2}\biggr).
\end{split}
\end{equation*}
Now substitute the geometric series
\begin{equation*}
\frac{x}{1+x^2}\biggl(\frac{x^2}{1+x^2}-\frac{1}{2}\biggr)
=-\frac{1}{2}\sum_{k=1}^{\infty}k\chi_{-4}(k)x^k,
\end{equation*}
and then swap the order of summation, to obtain
\begin{equation*}
\begin{split}
\int_{0}^{K}\frac{\cn^2 u\,\dn^2u}{1-\alpha\sn^4u}\,\log\sn u\,\d u
&=-\frac{\pi}{16}\log\biggl(\frac{\pi^4}{K^4}\cdot q\prod_{k=1}^{\infty}\frac{(1-q^{2k})^{16}}{(1-q^k)^8}\biggr)
\\ &\qquad
+\frac{\pi}{8}\biggl(\frac{\log q}{2}+2\sum_{k=1}^{\infty}k\chi_{-4}(k)\log(1-q^k)\biggr)
\\ &\qquad
-\frac{\pi}{4\sqrt{\alpha}}\sum_{k=1}^{\infty}k\chi_{-4}(k)\log\frac{1-q^k}{(1-q^{k/2})^2}.
\end{split}
\end{equation*}
Finally, by the $q$-series expansion for
$m(4/\sqrt{\alpha})$ \cite[Entry (2-9)]{LR}, and by
\cite[pg.~124, Entries 12.2 and~12.3]{Be3}, this becomes
\begin{equation*}
\begin{split}
\int_{0}^{K}\frac{\cn^2u\,\dn^2u}{1-\alpha\sn^4u}\log\sn u\,\d u
&=-\frac{\pi}{16}\log\alpha-\frac{\pi}{8}m\biggl(\frac{4}{\sqrt{\alpha}}\biggr)
\\ &\qquad
+\frac{\pi}{8\sqrt{\alpha}}\biggl(m\biggl(\frac{4}{\sqrt{\alpha}}\biggr)-2m\biggl(\frac{4}{\sqrt{\alpha'}}\biggr)\biggr),
\end{split}
\end{equation*}
where $\alpha$ has degree $2$ over~$\alpha'$.  By the second degree
modular equation of Ramanujan \cite[pg.~215]{Be3}, we know that
\begin{equation*}
\alpha'=\frac{4\sqrt{\alpha}}{(1+\sqrt{\alpha})^2}.
\end{equation*}
Since $0\le \alpha\le 1$, we can apply a functional equation of
Kurokawa and Ochiai \cite{KO}, to obtain
\begin{equation*}
2m\biggl(\frac{4}{\sqrt{\alpha'}}\biggr)
=2m\bigl(2(\alpha^{1/4}+\alpha^{-1/4})\bigr)
=m(4\sqrt{\alpha})+m\biggl(\frac{4}{\sqrt{\alpha}}\biggr).
\end{equation*}
This last observation completes the proof of \eqref{sn mahler
measure reduction}.  Formulas \eqref{cn mahler measure reduction}
and \eqref{dn mahler measure redution} can be proved with an
identical method, except that they require Fourier expansions for
$\log\cn u$ and $\log\dn u$, respectively, \cite[pg.~917]{GR}.
\end{proof}

\begin{proof}[Alternative proof of Theorem~\textup{\ref{th-F(2,3)}}]
If we let $p\mapsto v^2/2$ in~\eqref{G1c}, then
\begin{equation*}
G(1)=\frac{\pi}{6}\int_{0}^{1}\frac{\sqrt{(1-v^2)(1-\frac14v^2)}}{(1-\frac14v^4)}
\,\log\frac{v^6(1-\frac14v^2)}{4(1-v^2)}\,\d v.
\end{equation*}
Theorem \ref{th-F(2,3)} follows immediately from combining an
elementary result
\begin{equation*}
\frac{\pi}{4}=\int_{0}^{1}\frac{\sqrt{(1-v^2)(1-\alpha v^2)}}{(1-\alpha v^4)}\,\d v
\end{equation*}
(consider a Taylor series in $\alpha$), with all three formulas in
Proposition~\ref{conductor 24 mahler measure integral
proposition}, and the known identity $m(8)=4m(2)$ \cite{LR}.
\end{proof}

\section{Conductor 20}
\label{s-cond20}

In this section we prove Boyd's conjectures for elliptic curves of
conductor~$20$. Recall~\cite{Ono} that such curves are associated to
the modular form $\eta^2(q^2) \eta^2(q^{10})$, so it follows that
$L(E_{20},2)=F(1,5)$. The first step is to use Ramanujan's modular
equations to relate $L(E_{20},2)$ to an elementary integral. The
elementary integral can then be reduced to Mahler measures by
substituting doubly-periodic elliptic functions, or by using
hypergeometric functions. Define $S(x)$ as follows:
\begin{equation*}
S(x):=\int_{0}^{1}q^{(1+x)/4}\psi^2(q^x)\bigl(\psi^2(q)-5q\psi^2(q^5)\bigr)\,\log q\,\frac{\d q}{q}.
\end{equation*}
We begin by expressing $L(E_{20},2)$ in terms of $S(x)$.

\begin{lemma}
The following relation is true:
\begin{equation}
-4L(E_{20},2)=S(1)-S(5). \label{F(1,5) in terms of S}
\end{equation}
\end{lemma}

\begin{proof}
First notice that
\begin{equation*}
S(1)-S(5)
=\int_{0}^{1}q^{1/2}\bigl(\psi^2(q)-q\psi^2(q^5)\bigr)
\bigl(\psi^2(q)-5q\psi^2(q^5)\bigr)\,\log q\,\frac{\d q}{q}.
\end{equation*}
Ramanujan showed \cite[pg.~28]{BA} that
\begin{equation*}
\eta^2(q)\eta^2(q^5)=q^{1/2}\bigl(\psi^2(q)-q\psi^2(q^5)\bigr)\bigl(\psi^2(q)-5q\psi^2(q^5)\bigr),
\end{equation*}
which implies
\begin{equation*}
S(1)-S(5)=\int_{0}^{1}\eta^2(q)\, \eta^2(q^5)\,\log q\,\frac{\d
q}{q}=-4L(E_{20},2).
\end{equation*}
\vskip-\baselineskip
\end{proof}

Next we apply our trick to obtain a transformation for~$S(x)$.

\begin{proposition}
\label{prop5}
Suppose that $x>0$. Then
\begin{equation}\label{S(x) reduced}
S(x)=-\pi\int_{0}^{1}q^{x/2}\psi^4(-q^x)\log\biggl(5\frac{\varphi^2(q^{5})}{\varphi^2(q)}\biggr)\frac{\d q}{q}.
\end{equation}
\end{proposition}

\begin{proof}
Begin by setting $q=e^{-2\pi u}$, then
\begin{equation*}
S(x)=-(2\pi)^2\int_{0}^{\infty}u e^{-\pi xu/2}\psi^2(e^{-2\pi x u})
\bigl(e^{-\pi u/2}\psi^2(e^{-2\pi u})-5e^{-5\pi u/2}\psi^2(e^{-10\pi u})\bigr)\,\d u.
\end{equation*}
We use the following Lambert series expansion (which follows from
\cite[pg.~139, Example 4]{Be3}):
\begin{equation*}
e^{-\pi x u/2}\psi^2(e^{-2\pi x u})
=\sum_{n,k=1}^{\infty}\chi_{-4}(n)(e^{-\pi n k x u/2}-e^{-\pi n k x u}).
\end{equation*}
By the involution for the psi function and by \cite[pg.~114, Entry 8.1]{Be3}, we have
\begin{align}
e^{-\pi u/2}\psi^2(e^{-2\pi u})-5e^{-5\pi u/2}\psi^2(e^{-10\pi u})
&=\frac{1}{4u}\bigl(\varphi^2(-e^{-\pi/u})-\varphi^2(-e^{-\pi/(5u)})\bigr)
\notag\\
&=\frac{1}{u}\sum_{r,s=1}^{\infty}(-1)^s \chi_{-4}(r)(e^{-\pi r s/u}-e^{-\pi r s/(5u)}).
\notag
\end{align}
Therefore, the integral becomes
\begin{equation*}
\begin{split}
S(x)
=-(2\pi)^2\sum_{n,k,r,s\ge 1}(-1)^s\chi_{-4}(n r)
&\int_{0}^{\infty}(e^{-\pi n k x u/2}-e^{-\pi n k x u})
\\ &\qquad\times
(e^{-\pi r s/(u)}-e^{-\pi r s/(5 u)})\,\d u.
\end{split}
\end{equation*}
Use linearity and a $u$-substitution, to regroup the integral:
\begin{equation*}
\begin{split}
S(x)=-(2\pi)^2\sum_{n,k,r,s\ge 1}(-1)^s\chi_{-4}(n r)
&\int_{0}^{\infty}e^{-\pi n k x u}(2e^{-\pi r s/(2 u)}-2e^{-\pi r s/(10 u)}
\\ &\qquad
-e^{-\pi r s/u}+e^{-\pi r s/(5 u)})\, \d u.
\end{split}
\end{equation*}
Finally make the $u$-substitution $u\mapsto r u/k$. This
permutes the indices of summation inside the integral. We have
\begin{equation*}
\begin{split}
S(x)=-(2\pi)^2\sum_{n,k,r,s\ge 1}r\chi_{-4}(n r)\frac{(-1)^s}{k}
&\int_{0}^{\infty}e^{-\pi n r x u}(2e^{-\pi k s/(2 u)}-2e^{-\pi k s/(10 u)}
\\ &\qquad
-e^{-\pi k s/u}+e^{-\pi k s/(5 u)})\, \d u.
\end{split}
\end{equation*}
Simplify the $k$ and $s$ sums, then use the involution for the
Dedekind eta function and the product expansion
$\varphi(q)=\eta^5(q^2)/(\eta^2(q)\eta^2(q^4))$, to reduce the integral to
\begin{equation}\label{S(x) almost done}
S(x)=-2\pi^2\int_{0}^{\infty}\biggl(\sum_{n,r=1}^{\infty}r\chi_{-4}(r n)e^{-\pi r n x u}\biggr)
\log\biggl(5\frac{\varphi^2(e^{-10 \pi u})}{\varphi^2(e^{-2\pi u})}\biggr)\d u.
\end{equation}
Finally, the nested sum is easy to simplify.  By \cite[pg.~139, Example 3]{Be3},
\begin{equation*}
\sum_{n,r=1 }^{\infty}r\chi_{-4}(r n)q^{r n}
=\sum_{r=1}^{\infty}\frac{r\chi_{-4}(r)q^{r}}{1+q^{2r}}
=q\psi^4(-q^2).
\end{equation*}
Substituting this last result into \eqref{S(x) almost done}
completes the proof of~\eqref{S(x) reduced}.
\end{proof}

Now we can derive an elementary integral for $L(E_{20},2)$. In order
to accomplish the reduction, we need several additional modular
equations.

\begin{lemma}
We have
\begin{equation}\label{F(1,5) elementary integral}
L(E_{20},2)=-\frac{\pi}{20}\int_{0}^{1}\frac{(1-6t)\log(1+4t)}{\sqrt{t(1-t)(1+4t^2)}}\,\d
t.
\end{equation}
\end{lemma}

\begin{proof}
By formulas \eqref{S(x) reduced} and \eqref{F(1,5) in terms of
S}, we find that
\begin{equation*}
L(E_{20},2)=\frac{\pi}{4}\int_{0}^{1}\frac{q^{1/2}\psi^4(-q)-q^{5/2}\psi^4(-q^5)}{q}
\log\biggl(5\frac{\varphi^2(q^5)}{\varphi^2(q)}\biggr)\,\d q.
\end{equation*}
Now set $m=\varphi^2(q)/\varphi^2(q^5)$.  Then by \cite[pg.~26, formula (1.6.4)]{BA}, we obtain
\begin{equation*}
1-q^2\frac{\psi^4(-q^5)}{\psi^4(-q)}=\frac{8(3-m)}{(5-m)^2}.
\end{equation*}
Therefore, the integral becomes
\begin{equation*}
L(E_{20},2)=2\pi\int_{0}^{1}\frac{3-m}{(5-m)^2}
\log\biggl(\frac{5}{m}\biggr)\frac{q^{1/2}\psi^4(-q)}{q}\,\d q.
\end{equation*}
Now set $\alpha=1-\varphi^4(-q)/\varphi^4(q)$; then it is
known that
\begin{equation*}
\begin{split}
\frac{q^{1/2}\psi^4(-q)}{q}
&=\frac{1}{2\sqrt{4\alpha(1-\alpha)}}\,\frac{\d\alpha}{\d q}
\\
&=\frac{1}{8\sqrt{4\alpha(1-\alpha)(1-4\alpha(1-\alpha))}}\,
\frac{\d}{\d q}\bigl(4\alpha(1-\alpha)\bigr).
\end{split}
\end{equation*}
Finally, we have the following relation between $\alpha$ and~$m$:
\begin{equation*}
4\alpha(1-\alpha)=\frac{(m-1)(5-m)^5}{64 m^5}.
\end{equation*}
This relation between $\alpha$ and $m$ follows from
\cite[pg.~288, Entry 14]{Be3}.  Notice that the entry holds for $|q|<1$ by
the principle of analytic continuation. Eliminating $\alpha$, and
exercising caution about the square root, reduces the integral to
\begin{equation*}
\begin{split}
L(E_{20},2)
&=\frac{\pi}{4}\int_{0}^{1}\log\biggl(\frac{5}{m}\biggr)\frac{3-m}{m\sqrt{(5-m)(m-1)(5-2m+m^2)}}
\,\frac{\d m}{\d q}\,\d q
\\
&=\frac{\pi}{4}\int_{1}^{5}\log\biggl(\frac{5}{m}\biggr)\frac{3-m}{m\sqrt{(5-m)(m-1)(5-2m+m^2)}}\,\d m.
\end{split}
\end{equation*}
The change of variables from $q$ to $m$ is justified because $m$
ranges monotonically between $m=1$ and $m=5$ when $q\in[0,1]$.
Finally, the substitution $m\mapsto 5/(1+4t)$ completes the
proof of~\eqref{F(1,5) elementary integral}.
\end{proof}

Below we use two methods to reduce \eqref{F(1,5) elementary
integral} to Mahler measures.  The first method is to substitute
doubly-periodic elliptic functions into the integral. The main
draw-back to this method is that we first have to construct
non-standard elliptic functions.  The second method is to prove the
identity directly via hypergeometric manipulations. In both
approaches we investigate an integral which generalizes
\eqref{F(1,5) elementary integral}.  Notice that
\begin{equation}
J(y):=\frac1{2\pi}\int_0^1\frac{(2-y+3yt)\log(1+yt)}{\sqrt{t(1-t)(4+(4-y)yt+y^2t^2)}}\,\d t
\label{J(y)}
\end{equation}
reduces to the integral in~\eqref{F(1,5) elementary integral} when
$y=4$.

\subsection{The elliptic reduction}
Throughout this subsection we assume that $k>4/3$.  Notice that
when $y=2k/(k-1)$, we have
\begin{equation}\label{J(2k/(k-1)) integral definition}
J\biggl(\frac{2k}{k-1}\biggr)
=-\frac1{\pi}\int_0^1\frac{(1-3kt)\log\biggl(1-\dfrac{2kt}{1-k}\biggr)}
{\sqrt{4t((1-k)^2-t(1-kt)^2)}}\,\d t.
\end{equation}
The overarching goal of the following discussion is to obtain
formula \eqref{I(k) general integral formula}.  To prove that
identity, it is necessary to use Fourier series expansions for
elliptic functions which parameterize the curve
\begin{equation*}
F_k: y^2=4x\bigl((1-k)^2-x(1-k x)^2\bigr).
\end{equation*}
Since we (regrettably) could not find such formulas in the
literature, we first prove Proposition \ref{proposition on values of
w(x)} and Lemma~\ref{w(x) fourier series lemma}.

Notice that $F_k$ is a genus-one curve, with non-zero discriminant
when $k>4/3$. Therefore, $F_k$ can be parameterized by
doubly-periodic functions. Suppose that $w(x)$ satisfies the
differential equation:
\begin{equation}\label{w(x) differential equation}
\bigl(w'(x)\bigr)^2=4w(x)\bigl((1-k)^2-w(x)(1-k w(x))^2\bigr).
\end{equation}
In order to explicitly identify $w(x)$ we can map $F_k$ to
$Y^2=4X^3-g_2 X-g_3$.  It follows easily that
\begin{equation}\label{w in terms of weierstrass p}
w(x)=\frac{3(1-k)^2}{1+3\wp(x)},
\end{equation}
where $\wp(x):=\wp(x,\{g_2,g_3\})$ is the Weierstrass function,
and
\begin{align*}
g_2&=-\frac{4}{3}(6k^3-12k^2+6k-1),
\displaybreak[2]\\
g_3&=\frac{4}{27}(2-6k+3k^2)(1-6k+12k^2-18k^3+9k^4).
\end{align*}
This identification is quite useful for computational purposes.

\begin{proposition}\label{proposition on values of w(x)}
Let $2K$ and $2K'$ denote the real and purely imaginary periods of~$w(x)$.
Then $w(x)$ has the following values:
\begin{equation*}
    \begin{tabular}{|c|c|c|c|p{6 in}|}
        \hline
        $x$ &  $w(x)$ & order of zero/pole &  residue\\
        \hline
        $0$ & $0$ & $2$&  $-$\\
        $K$ & $1$ &   $-$& $-$\\
        $K'$ & $1$ &  $-$ & $-$\\
        $K+K'$ & $0$ &  $2$ & $-$\\
        $K'/3$ & $(1-k)/(2k)$ & $-$ & $-$\\
        $2K'/3$ & $\infty$ &  $1$ & $i/(2k)$\\
        $4K'/3$ & $\infty$ & $1$ & $-i/(2k)$\\
        $5K'/3$ & $(1-k)/(2k)$ &  $-$ & $-$\\
        $K+K'/3$ & $\infty$ & $1$ & $-i/(2k)$\\
        $K+2K'/3$ & $(1-k)/(2k)$ &  $-$ & $-$\\
        $K+4K'/3$ & $(1-k)/(2k)$ & $-$ & $-$\\
        $K+5K'/3$ & $\infty$ &  $1$ & $i/(2k)$\\
        \hline
    \end{tabular}
\end{equation*}
\end{proposition}

\begin{proof}
It is well known that $\wp(x)$ has a second-order pole at $x=0$,
so $w(x)$ has a second-order zero at that point. Since $\wp(x)$ is
even, $w(2a K+2b K'-x)=w(x)$ for all $(a,b)\in\mathbb{Z}^2$.
Therefore we only need to evaluate $w(x)$ for $x\in\{K, K', K+K',
K'/3,2K'/3, K+4K'/3, K+5K'/3\}$.

We will require additional properties of the Weierstrass
$\wp$-function.  If $4X^3-g_2 X-g_3=4(x-r_1)(x-r_2)(x-r_3)$, then
the half-periods of $\wp(x)$ are given by
\begin{equation*}
\omega=\int_{\infty}^{r_1}\frac{1}{\sqrt{4y^3-g_2 y-g_3}}\,\d y,
\qquad
\omega'=\int_{r_1}^{r_2}\frac{1}{\sqrt{4y^3-g_2 y-g_3}}\,\d y.
\end{equation*}
Select $r_1=(1-k)^2-\frac{1}{3}$ to be the real zero of $4X^3-g_2
X-g_3=0$, and $r_2$ to be the imaginary zero which lies in the
upper half plane.  Now set $K:=\omega$, and $K':=2\omega-\omega'$.
While it is possible to show $\Re K'=0$, we will not pursue that
calculation here. It follows that $K+K'$ is a period of $w(x)$, so
we have the following identities:
\begin{align*}
w(K+K')&=w(0)=0,\\
w(K)&=w(2K+K')=w(K'),\\
w(K'/3)&=w(K+4K'/3),\\
w(2K'/3)&=w(K+5K'/3).
\end{align*}
Since $\wp(\omega)=r_1=(1-k)^2-\frac{1}{3}$, we can use \eqref{w in
terms of weierstrass p} to conclude that $w(K)=1$. The values of
$w(K'/3)$ and $w(2K'/3)$ can be verified from a polynomial relation
between $w(x)$ and $w(3x)$, which follows from the Weierstrass
addition formula. Now we calculate the values of the residues. Since
$w(2K'/3)=\infty$, it follows that
$\wp\bigl(\frac{4\omega'-2\omega}{3}\bigr)=-\frac13$, thus
$\bigl(\wp'\bigl(\frac{4\omega'-2\omega}{3}\bigr)\bigr)^2=-4k^2(1-k)^4$.
Extracting a square root we obtain
$\wp'\bigl(\frac{4\omega'-2\omega}{3}\bigr)=-2i k(1-k)^2$. The
choice of square root can be justified by checking the formula
numerically at $k=2$, and then appealing to the fact that $\omega$,
$\omega'$, and $\wp'\bigl(\frac{4\omega'-2\omega}{3}\bigr)$ are
analytic functions of $k$ for $k>4/3$. Finally, by formula \eqref{w
in terms of weierstrass p}
\begin{equation*}
\res_{x=2K'/3}w(x)=\frac{(1-k)^2}{\wp'\bigl(\frac{4\omega'-2\omega}{3}\bigr)}=\frac{i}{2k}.
\end{equation*}
The other residues can be verified in a similar fashion.
\end{proof}

Notice that we can integrate \eqref{w(x) differential equation},
and use $w(K)=1$, to obtain a second formula for~$K$:
\begin{equation*}
K=\int_{0}^{1}\frac{\d t}{\sqrt{4t((1-k)^2-t(1-k t)^2)}}.
\end{equation*}
We need two Fourier series expansions to finish the elliptic
reduction.

\begin{lemma}\label{w(x) fourier series lemma}
Suppose that $x>0$.  Then
\begin{equation}\label{w(x) Fourier series}
w(x)=\frac{2\pi}{k K}\sum_{n=1}^{\infty}\frac{(-1)^{n+1}q^n}{1+(-1)^n q^n+q^{2n}}
\,\sin^2\frac{\pi n x}{2K}
\end{equation}
and
\begin{equation}\label{log(1+4w(x)) Fourier series}
\log\biggl(1-\frac{2k}{1-k}w(x)\biggr)
=8\sum_{\substack{n=1\\n \text{ odd}}}^{\infty}
\frac{1}{n}\frac{q^{n}-q^{2n}}{1+q^{3n}}\,\sin^2\frac{\pi n x}{2K},
\end{equation}
where $q=e^{2\pi i K'/6K}$.  An alternative formula for $q$ is
given by
\begin{equation}\label{w(x) fourier series value of q}
q=\exp\biggl(-\frac{2\pi}{\sqrt{3}}\frac{{}_2F_1(\frac{1}{3},\frac{2}{3};1;1-\alpha)}
{{}_2F_1(\frac{1}{3},\frac{2}{3};1;\alpha)}\biggr),
\end{equation}
where
\begin{equation*}
\alpha=\frac{27 p (1 + p)^4}{2 (1 + 4 p + p^2)^3},
\qquad
p=\frac{-1+\sqrt{(3k-1)/(k-1)}}{2}.
\end{equation*}
\end{lemma}

\begin{proof}
The proof of \eqref{w(x) Fourier series} is an
exercise in the theory of elliptic functions.  The poles of $w(x)$
inside the fundamental parallelogram are $2K'/3$, $4K'/3$,
$K+K'/3$, and $K+5K'/3$.  The function has residues $i/(2k)$,
$-i/(2k)$, $-i/(2k)$, and $i/(2k)$ at each of the poles. We also know
that $w(0)=w'(0)=0$.  By \cite[formula (27)]{BF}, we have
\begin{equation*}
\begin{split}
w(x)=\frac{3i}{2k}\sum_{\substack{m,n=-\infty\\
(m,n)\ne(0,0)}}^{\infty}\chi_{-3}(n)&\biggl(\frac{1}{3x-(6mK+2nK')}
\\ &\quad
+\frac{1}{(6m K+2nK')}+\frac{3x}{(6m K+2nK')^2}\biggr)
\\
-\frac{3i}{2k}\sum_{\substack{m,n=-\infty\\
(m,n)\ne(0,0)}}^{\infty}\chi(n)&\biggl(\frac{1}{3x-((6m+3)K+n K')}
\\ &\quad
+\frac{1}{((6m+3) K+n K')}+\frac{3x}{((6m+3) K+nK')^2}\biggr),
\end{split}
\end{equation*}
where $\chi(n)$ is the Legendre symbol mod $6$. It is a lengthy
exercise to reduce this last expression to \eqref{w(x) Fourier
series}. The fastest (if least rigorous) method for finishing the
calculation, is to differentiate the entire expression twice, and
then substitute the following Fourier series:
\begin{equation*}
\sum_{n=-\infty}^{\infty}\frac{1}{(x+\tau+n)^3}+\frac{1}{(-x+\tau+n)^3}
=i(2\pi)^3\sum_{n=1}^{\infty}n^2e^{2\pi i n\tau}\cos(2\pi n x),
\end{equation*}
which holds for $\Im(\tau)>0$.  Thus one obtains a formula for
$w''(x)$, which can be integrated to recover \eqref{w(x) Fourier
series}.

Proposition~\ref{proposition on values of w(x)} shows that
$(1-k)/(2k)=w(K'/3)$.  It follows that
\begin{equation*}
1-\frac{2k}{1-k}w(x)=1-\frac{w(x)}{w(K'/3)}.
\end{equation*}
This function has simple zeros at $K'/3$, $5K'/3$, $K+2K'/3$ and
$K+4K'/3$, and simple poles at $2K'/3$, $4K'/3$, $K+K'/3$ and
$K+5K'/3$. Since any two elliptic functions with the same zeros
and poles are constant multiples, it is easy to obtain an infinite
product. We have
\begin{equation}\label{1-2k/(1-k)*w(x) theta product}
1-\frac{2k}{1-k}w(x)
=C\frac{\theta(x,K'/3)\theta(x,5K'/3)\theta(x,K+2K'/3)\theta(x,K+4K'/3)}
{\theta(x,2K'/3)\theta(x,4K'/3)\theta(x,K+K'/3)\theta(x,K+5K'/3)},
\end{equation}
where
\begin{equation*}
\theta(x,\rho)=\bigl(1-e^{2\pi i(x-\rho)/(2K)}\bigr)
\prod_{n=1}^{\infty}\bigl(1-e^{2\pi i(x-\rho+2n K')/(2K)}\bigr)
\bigl(1-e^{2\pi i(-x+\rho+2n K')/(2K)}\bigr).
\end{equation*}
The right-hand side of \eqref{1-2k/(1-k)*w(x) theta product} is
doubly periodic because $\theta(x,\rho)$ has period $2K$,  and
satisfies the quasi-periodicity relation
$$
\theta(x+2K',\rho)=-e^{2\pi i(\rho-x)/(2K)}\theta(x,\rho).
$$
The right-hand side also has the correct zeros and poles, since
$\theta(x,\rho)$ vanishes at~$\rho$. The constant $C$ can be
determined by using the fact that $w(0)=0$. Finally,
\eqref{log(1+4w(x)) Fourier series} follows from taking logarithms
of \eqref{1-2k/(1-k)*w(x) theta product}, and then using the
Taylor series for the logarithm.

We conclude the proof by simplifying the expression for~$q$. Since
$w(K)=1$, we can use \eqref{log(1+4w(x)) Fourier series} to obtain
\begin{align*}
\log\frac{3k-1}{k-1}
&=8\sum_{n=0}^{\infty}\frac{1}{2n+1}\frac{q^{2n+1}-q^{4n+2}}{1+q^{6n+3}}
\\
&=4\log\prod_{n=1}^{\infty}\frac{(1-q^{2n})^5(1-q^{3n})^2(1-q^{12n})^2}{(1-q^{n})^2(1-q^{4n})^2(1-q^{6n})^5}
\\
&=2\log\frac{\varphi^2(q)}{\varphi^2(q^3)}.
\end{align*}
If we let $1+2p=\varphi^2(q)/\varphi^2(q^3)$, then it follows
easily that
$$
p=\frac{-1+\sqrt{(3k-1)/(k-1)}}{2}.
$$
Finally, formula \eqref{w(x) fourier series value of q} is a consequence of
standard inversion formulas in the theory of signature~$3$.
\end{proof}

\begin{theorem}
\label{elliptic-cond-20}
Suppose that $k\ge4/3$, and let $p=\frac12(-1+\sqrt{(3k-1)/(k-1)})$.  The following formula is
true:
\begin{equation}\label{I(k) general integral formula}
J\biggl(\frac{2k}{k-1}\biggr)
=2g\biggl(\frac{2(1+p)^2}{p}\biggr)-g\biggl(\frac{4(1+p)}{p^2}\biggr).
\end{equation}
\end{theorem}

\begin{proof}
First assume that $k>4/3$.  If we set $t=w(x)$, then
\eqref{J(2k/(k-1)) integral definition} becomes
\begin{equation*}
J\biggl(\frac{2k}{k-1}\biggr)
=-\frac{1}{\pi}\int_{0}^{K}(1-3 kw(x))\log\biggl(1-\frac{2k}{1-k}w(x)\biggr)\d x.
\end{equation*}
Substituting \eqref{log(1+4w(x)) Fourier series} and \eqref{w(x)
Fourier series} reduces the integral to
\begin{equation*}
\begin{split}
J\biggl(\frac{2k}{k-1}\biggr)
&=-\biggl(\frac{4K}{\pi}+12\sum_{j=1}^{\infty}\frac{(-q)^{j}}{1+(-q)^j+q^{2j}}\biggr)
\sum_{\substack{n=1\\\text{$n$ odd}}}^{\infty}\frac{1}{n}\,\frac{q^n-q^{2n}}{1+q^{3n}}
\\ &\qquad
+6\sum_{\substack{n=1\\\text{$n$ odd}}}^{\infty}\frac{1}{n}\,\frac{q^{2n}(1-q^{2n})}{(1+q^{3n})^2}.
\end{split}
\end{equation*}
Now substitute
\begin{equation*}
3\frac{q^2(1-q^2)}{(1+q^3)^2}=-\frac{q(1-q)}{(1+q^3)}-\sum_{j=1}^{\infty}(-1)^j j\chi_{-3}(j)q^j,
\end{equation*}
to obtain
\begin{align*}
J\biggl(\frac{2k}{k-1}\biggr)
&=-\biggl(\frac{4K}{\pi}+2a(-q)\biggr)\sum_{\substack{n=1\\\text{$n$ odd}}}^{\infty}
\frac{1}{n}\,\frac{q^n-q^{2n}}{1+q^{3n}}
\\ &\qquad
-\sum_{j=1}^{\infty}(-1)^jj\chi_{-3}(j)\log\frac{1+q^j}{1-q^j}.
\end{align*}
Notice that we have used a Lambert series for $a(-q)$, which
follows from \cite[pg.~100, Theorem 2.12]{Be5}.  Now we claim that
$2K=-\pi a(-q)$.  If one substitutes the hypergeometric
representation for $a(-q)$, then this statement is equivalent to
Lemma~\ref{L1} in the next subsection.  It is also possible to
prove the equality directly by setting $x=K$ in \eqref{w(x)
Fourier series} and \eqref{log(1+4w(x)) Fourier series}, and then
performing $q$-series manipulations. The $q$-series for
$J(2k/(k-1))$ reduces to
\begin{equation*}
J\biggl(\frac{2k}{k-1}\biggr)
=-\sum_{j=1}^{\infty}(-1)^j j\chi_{-3}(j)\log\frac{1+q^j}{1-q^j}.
\end{equation*}
We can now substitute Stienstra's $q$-series for $g(k)$ \cite{St}:
applying formula (2-11) in~\cite{LR} completes the proof
of~\eqref{I(k) general integral formula} if $k>4/3$. Finally,
notice that both sides of \eqref{I(k) general integral formula}
are continuous at $k=4/3$, hence the formula remains true for the
boundary value as well.
\end{proof}

\subsection{The hypergeometric reduction}

In this part, we show the coincidence of the derivatives of $J(y)$
and $g(y)$ on the interval $2<y<8$ and conclude with the identity
$J(y)=g(y)$ for $2\le y\le 8$ by appealing to the equality at
$y=8$ deduced in Theorem~\ref{elliptic-cond-20}.

\begin{lemma}
\label{L1}
For $2\le y<8$, we have
\begin{equation}
\frac1{2\pi}\int_0^1\frac{\d t}{\sqrt{t(1-t)(4+(4-y)yt+y^2t^2)}}
=\frac1{y+4}\,{}_2F_1\biggl(\begin{matrix} \frac13, \, \frac23 \\
1 \end{matrix}\biggm| \frac{27y^2}{(y+4)^3} \biggr).
\label{T5}
\end{equation}
\end{lemma}

\begin{proof}
We apply the transformation \cite[p.~112, Theorem~5.6]{Be5},
\begin{equation}
\frac1{1+p+p^2}\,{}_2F_1\biggl(\begin{matrix} \frac13, \, \frac23 \\
1 \end{matrix}\biggm| \frac{27p^2(1+p)^2}{4(1+p+p^2)^3} \biggr)
=\frac1{\sqrt{1+2p}}\,{}_2F_1\biggl(\begin{matrix} \frac12, \, \frac12 \\
1 \end{matrix}\biggm| \frac{p^3(2+p)}{1+2p} \biggr)
\label{T5a}
\end{equation}
with the choice $p=(\sqrt{1+y}-1)/2$ (ranging in $(\sqrt3-1)/2\le p<1$), so that $y=4p(1+p)$
and the left-hand side in~\eqref{T5a} assumes the form
\begin{equation}
\frac1{1+p+p^2}\,{}_2F_1\biggl(\begin{matrix} \frac13, \, \frac23 \\
1 \end{matrix}\biggm| \frac{27p^2(1+p)^2}{4(1+p+p^2)^3} \biggr)
=\frac4{y+4}\,{}_2F_1\biggl(\begin{matrix} \frac13, \, \frac23 \\
1 \end{matrix}\biggm| \frac{27y^2}{(y+4)^3} \biggr).
\label{T5b}
\end{equation}
On the other hand, the substitution $y=4p(1+p)$ and the change of variable
$$
t=\frac{1-u}{1+2pu-p^3(2+p)u^2}
$$
in the original integral results in
\begin{align}
\frac1\pi\int_0^1\frac{\d t}{\sqrt{t(1-t)(4+(4-y)yt+y^2t^2)}}
&=\frac1{2\pi}\int_0^1\frac{\d u}{\sqrt{u(1-u)(1+2p-p^3(2+p)u)}}
\notag\\
&=\frac1{2\sqrt{1+2p}}\,{}_2F_1\biggl(\begin{matrix} \frac12, \, \frac12 \\
1 \end{matrix}\biggm| \frac{p^3(2+p)}{1+2p} \biggr),
\label{T5c}
\end{align}
where on the last step we apply the Euler--Pochhammer integral representation
of the hypergeometric series \cite[equation (1.6.6)]{Slater}. Combining \eqref{T5a}--\eqref{T5c}
we arrive at the desired claim~\eqref{T5}.
\end{proof}

Note the range of the argument of the hypergeometric series in~\eqref{T5}:
$$
\frac12\le\frac{27y^2}{(y+4)^3}<1
\quad\text{for $2\le y<8$}.
$$

\begin{lemma}
\label{L2}
For $2<y<8$,
\begin{equation*}
\frac{\d J}{\d y}
=\frac1{y+4}\,{}_2F_1\biggl(\begin{matrix} \frac13, \, \frac23 \\
1 \end{matrix}\biggm| \frac{27y^2}{(y+4)^3} \biggr).
\end{equation*}
\end{lemma}

\begin{proof}
Note that for real values of $y$ in the interval $2\le y\le8$
we have
$$
\biggl|\sqrt t\biggl(1-\frac y2(1-t)\biggr)\biggr|\le1
\quad\text{for $0\le t\le1$},
$$
so that the real-valued function
\begin{equation*}
v(t)=v(t;y)=2\arcsin\biggl(\sqrt t\biggl(1-\frac y2(1-t)\biggr)\biggr)
\end{equation*}
is well defined on the interval $0<t<1$. Because
\begin{equation}
\frac{\partial v}{\partial t}=\frac{2-y+3yt}{\sqrt{t(1-t)(4+(4-y)yt+y^2t^2)}},
\label{T2}
\end{equation}
we can write the integral~\eqref{J(y)} as
\begin{equation}
 J(y)=\frac1{2\pi}\int_0^1\log(1+yt)\frac{\partial v}{\partial t}\,\d t. \label{T3}
\end{equation}
Denote $u(t)=u(t;y)=\log(1+yt)$ and use, besides~\eqref{T2},
$$
\frac{\partial u}{\partial t}=\frac y{1+yt}, \quad
\frac{\partial u}{\partial y}=\frac t{1+yt}, \quad\text{and}\quad
\frac{\partial v}{\partial y}=-\frac{2\sqrt{t(1-t)}}{\sqrt{4+(4-y)yt+y^2t^2}}.
$$
It follows from~\eqref{T3} that
\begin{align*}
\frac{\d}{\d y} J(y) &=\frac1{2\pi}\int_0^1\frac\partial{\partial
y}\biggl(u\frac{\partial v}{\partial t}\biggr)\,\d t
=\frac1{2\pi}\int_0^1\biggl(\frac{\partial u}{\partial
y}\,\frac{\partial v}{\partial t} +u\frac{\partial^2v}{\partial
y\,\partial t}\biggr)\,\d t
\\
&=\frac1{2\pi}\int_0^1\frac{\partial u}{\partial y}\,\frac{\partial v}{\partial t}\,\d t
+\frac1{2\pi}\int_0^1u\,\d\biggl(\frac{\partial v}{\partial y}\biggr)
\intertext{(integrating the second integral by parts)}
&=\frac1{2\pi}\int_0^1\frac{\partial u}{\partial y}\,\frac{\partial v}{\partial t}\,\d t
+\frac1{2\pi}\,u\frac{\partial v}{\partial y}\bigg|_{t=0}^{t=1}
-\frac1{2\pi}\int_0^1\frac{\partial v}{\partial y}\,\frac{\partial u}{\partial t}\,\d t
\\
&=\frac1{2\pi}\int_0^1\biggl(\frac{\partial u}{\partial y}\,\frac{\partial v}{\partial t}
-\frac{\partial v}{\partial y}\,\frac{\partial u}{\partial t}\biggr)\,\d t
\\
&=\frac1{2\pi}\int_0^1\frac{t(2+y+yt)}{(1+yt)\sqrt{t(1-t)(4+(4-y)yt+y^2t^2)}}\,\d t
\\
&=\frac1{2\pi}\int_0^1\frac{t\,\d t}{\sqrt{t(1-t)(4+(4-y)yt+y^2t^2)}}
\\ &\qquad
+\frac1{2\pi}\int_0^1\frac{(1+y)t\,\d t}{(1+yt)\sqrt{t(1-t)(4+(4-y)yt+y^2t^2)}}
\intertext{(applying the change $t\mapsto(1-t)/(1+yt)$ in the second integral)}
&=\frac1{2\pi}\int_0^1\frac{t\,\d t}{\sqrt{t(1-t)(4+(4-y)yt+y^2t^2)}}
\\ &\qquad
+\frac1{2\pi}\int_0^1\frac{(1-t)\,\d t}{\sqrt{t(1-t)(4+(4-y)yt+y^2t^2)}}
\\
&=\frac1{2\pi}\int_0^1\frac{\d t}{\sqrt{t(1-t)(4+(4-y)yt+y^2t^2)}}.
\end{align*}
It remains to apply Lemma~\ref{L1} to the resulting integral.
\end{proof}

\begin{theorem}
\label{th1}
For $2\le y\le 8$, the equality
\begin{equation}
J(y)=g(y)
\label{J=g}
\end{equation}
holds.
\end{theorem}

\begin{proof}
For $2<y<8$, the hypergeometric evaluation~\eqref{hyper-g} of $g(y)$ can be stated
in the form
$$
g(y)=\frac13f\biggl(\frac{y^2}{(y+4)^3}\biggr)+\frac43\Re f\biggl(\frac{y}{(y-2)^3}\biggr)
$$
where the function
\begin{equation}
f(z)=-\frac{\log z}3-2z\,{}_4F_3\biggl(\begin{matrix} \frac43, \, \frac53, \, 1, \, 1 \\
2, \, 2, \, 2 \end{matrix}\biggm| 27z \biggr)
\label{f-notation}
\end{equation}
satisfies the equation
$$
\frac{\d f}{\d z}=-\frac1{3z}\,{}_2F_1\biggl(\begin{matrix} \frac13, \, \frac23 \\
1 \end{matrix}\biggm| 27z \biggr).
$$
Therefore,
$$
\frac{\d g}{\d y}=\frac{y-8}{9y(y+4)}\,{}_2F_1\biggl(\begin{matrix} \frac13, \, \frac23 \\
1 \end{matrix}\biggm| \frac{27y^2}{(y+4)^3} \biggr)
+\frac{8(y+1)}{9y(y-2)}\,\Re{}_2F_1\biggl(\begin{matrix} \frac13, \, \frac23 \\
1 \end{matrix}\biggm| \frac{27y}{(y-2)^3} \biggr),
$$
and application of the cubic transformation
$$
\Re{}_2F_1\biggl(\begin{matrix} \frac13, \, \frac23 \\
1 \end{matrix}\biggm| \frac{27y}{(y-2)^3} \biggr)
=\frac{y-2}{y+4}\,{}_2F_1\biggl(\begin{matrix} \frac13, \, \frac23 \\
1 \end{matrix}\biggm| \frac{27y^2}{(y+4)^3} \biggr)
$$
result in
\begin{equation}
\frac{\d g}{\d y}=\frac1{y+4}\,{}_2F_1\biggl(\begin{matrix} \frac13, \, \frac23 \\
1 \end{matrix}\biggm| \frac{27y^2}{(y+4)^3} \biggr).
\label{dg/dy}
\end{equation}
Comparing this evaluation with the one from Lemma~\ref{L2} we
conclude that $g(y)$ and $ J(y)$ differ on the interval $2<y<8$ by
a constant; because both $g(y)$ and $ J(y)$ are continuous at the
end-points, the relation $ J(y)-g(y)=C$, a real constant, is true
for $2\le y\le8$. To determine the constant, take $y=8$ and apply
Theorem~\ref{elliptic-cond-20} with the choice $p=1$; it follows
that
$$
g(8)+C=J(8)=2g(8)-g(8),
$$
hence $C=0$.
\end{proof}

\begin{remark}
The derivative~\eqref{dg/dy} can be alternatively obtained
by differentiating Stienstra's $q$-series for $g(y)$ \cite[Example \#6]{St},
\cite[formula~(2-11)]{LR}.
\end{remark}

\subsection{Culmination}

We conclude this section by listing the major consequences
of Theorems~\ref{elliptic-cond-20} and~\ref{th1}.

\begin{theorem}
\label{cor0}
The following formulas are true:
\begin{align}
\frac{10}{\pi^2}L(E_{20},2) &=2g(4+2\sqrt{5})-g(8+4\sqrt{5})
\label{F(1,5) in terms of algebraic g}\\
&=g(4)\label{F(1,5) final elliptic reduction}\\
&=\frac34n(\sqrt[3]{32}).\label{F(1,5) in terms of algebraic n}
\end{align}
\end{theorem}

\begin{proof}
Equation \eqref{F(1,5) in terms of algebraic g}
follows from setting $k=2$ in \eqref{I(k) general integral
formula} and then comparing it to \eqref{F(1,5) elementary
integral}, while \eqref{F(1,5) final elliptic reduction}
follows from taking $y=4$ in~\eqref{J=g}.

Using the hypergeometric evaluations \eqref{hyper-g} and \eqref{hyper-n} in the form
$$
g(k)=\frac13f\biggl(\frac{k^2}{(k+4)^3}\biggr)+\frac43f\biggl(\frac{k}{(k-2)^3}\biggr),
\qquad
n(k)=f\biggl(\frac1{k^3}\biggr),
$$
whenever the arguments lie between 0 and $1/27$,
with the hypergeometric function $f(z)$ defined in~\eqref{f-notation},
we have
\begin{align*}
&
2g(4+2\sqrt{5})-g(8+4\sqrt{5})
\\ &\quad
=\frac23f\biggl(\frac{4}{(7-\sqrt5)^3}\biggr)+\frac83f\biggl(\frac1{32}\biggr)
-\frac13f\biggl(\frac1{32}\biggr)-\frac43f\biggl(\frac{4}{(7+\sqrt5)^3}\biggr)
\\ &\quad
=\frac73n(\sqrt[3]{32})+\frac23n\biggl(\frac{7-\sqrt5}{\sqrt[3]{4}}\biggr)
-\frac43n\biggl(\frac{7+\sqrt5}{\sqrt[3]{4}}\biggr).
\end{align*}
Finally, Bertin's ``exotic'' relation \cite[Theorem~6]{GuR}
\begin{equation}
16n\biggl(\frac{7+\sqrt5}{\sqrt[3]{4}}\biggr)
-8n\biggl(\frac{7-\sqrt5}{\sqrt[3]{4}}\biggr)
=19n(\sqrt[3]{32})
\label{exotic}
\end{equation}
reduces the latter sum to $\frac34n(\sqrt[3]{32})$.
\end{proof}

\begin{corollary}
\label{cor00}
The following Boyd's conjectural evaluations are true:
\begin{align*}
n(\sqrt[3]{2})=\frac{25}{6\pi^2}L(E_{20},2),
\qquad g(-2)=\frac{15}{\pi^2}L(E_{20},2).
\end{align*}
\end{corollary}

\begin{proof}
These readily follow from \cite[formula (2-26)]{LR},
\begin{align*}
3g(-2)=n(2^{1/3})+4n(2^{5/3}),
\qquad
3g(4)=4n(2^{1/3})+n(2^{5/3}),
\end{align*}
and Theorem~\ref{cor0}.
\end{proof}

\begin{theorem}
\label{cor1}
For $(\sqrt3-1)/2\le p\le 1$, the Mahler measure $g(\,\cdot\,)$ satisfies the functional equation
\begin{equation*}
g\bigl(4p(1+p)\bigr)+g\biggl(\frac{4(1+p)}{p^2}\biggr)=2g\biggl(\frac{2(1+p)^2}p\biggr).
\end{equation*}
\end{theorem}

\begin{proof}
The result follows by comparing the two different evaluations
obtained in Theorems~\ref{elliptic-cond-20} and~\ref{th1}.
\end{proof}

\begin{remark}
In view of the proof of Theorem~\ref{cor0}, our
Theorem~\ref{cor1} may be thought of as a generalization
of Bertin's ``exotic'' relation~\eqref{exotic}.
\end{remark}

\begin{theorem}
\label{cor2}
We have
\begin{equation*}
L(E_{36},2)
=-\frac{2\pi^2\log2}{27}+\frac{\Gamma^3(\tfrac13)}{3\cdot2^{7/3}}
\,{}_3F_2\biggl(\begin{matrix} \frac13, \, \frac13, \, 1 \\
\frac56, \, \frac43 \end{matrix}\biggm| -\frac18 \biggr)
+\frac{\Gamma^3(\tfrac23)}{2^{11/3}}
\,{}_3F_2\biggl(\begin{matrix} \frac23, \, \frac23, \, 1 \\
\frac76, \, \frac53 \end{matrix}\biggm| -\frac18 \biggr).
\end{equation*}
\end{theorem}

\begin{proof}
Rodriguez-Villegas \cite{RV} showed that
$$
L(E_{36},2)=\frac{2\pi^2}9g(2).
$$
Making the change $t^3=u$ in the integral
$$
L(E_{36},2)=\frac{2\pi^2}9 J(2) =\frac\pi3\int_0^1\frac{\sqrt
t\log(1+2t)}{\sqrt{1-t^3}}\,\d t
$$
and writing the interior logarithm as hypergeometric series,
we arrive at the claim. Note that both \texttt{Maple} and \texttt{Mathematica}
produce the evaluation without human assistance.
\end{proof}

\section{Concluding remarks}
\label{sec5}

We conclude by mentioning the fact that this paper settles the
conjectures of Bloch and Grayson for elliptic curves of conductor
$20$ \cite{BG}.  Since there is a very simple method to translate
Mahler measures into elliptic dilogarithms \cite{GuR}, the main
results in this paper are equivalent to relations between $L(E,2)$
and values of the elliptic dilogarithm. For instance, given the
conductor $20$ elliptic curve $E: y^2=4x^3-432 x+1188$, and the
torsion point $P=(-6,54)$, by \cite[~Theorem 4]{GuR} we have
\begin{equation*}
D^{E}(2P)=\frac{2\pi}{9}n(2^{5/3})=\frac{80}{27\pi}L(E,2).
\end{equation*}
The equality to $L(E,2)$ follows immediately from~\eqref{F(1,5) in
terms of algebraic n}.  Although Bloch and Grayson did not examine
any conductor $24$ curves, we can prove similar relations for
those cases, by combining Theorem~\ref{cor0} with \cite[~Theorem 4]{GuR}.

Finally, there are many additional problems which need to be
addressed. The most obvious direction is to try to prove more of
Boyd's conjectures. There are still hundreds of outstanding
conjectures in Boyd's tables \cite{Bo1}.  It would also be
interesting to understand what overlap (if any) exists between our
techniques, and those of Brunault \cite{Br} and Mellit~\cite{Me}.
They proved Boyd's conjectures for elliptic curves of conductors
$11$ and $14$ by using Beilinson's theorem. Rodriguez-Villegas was
the first to advocate this $K$-theoretic approach~\cite{RV}; he originally
suggested that the conductor~$24$ cases could be proved with
Beilinson's theorem.

It would also be interesting to reduce more values of $F(b,c)$ to
hypergeometric functions. An easy corollary to the $L(E_{20},2)$
formula of Theorem~\ref{cor0}, is a formula for $F(5,9)$. By
\cite{Rgsumbit2} we know that $9 F(5,9)=45 F(1,1)-50 F(1,5)$,
hence we obtain
\begin{equation}
\frac{18}{5\pi^2}F(5,9)=g(-4)-2g(4).
\end{equation}
Notice that Lemma \ref{Lemma on special values of H} and
Proposition \ref{Theorem on H(x) in elementary integrals} reduce
$F(3,7)$, $F(6,7)$, and $F(3/2,7)$ to complicated elementary
integrals. We expect these lattice sums to also equal values of
hypergeometric functions, although we currently see no way to
prove it.

\begin{acknowledgements}
The authors would like to thank Bruce Berndt, David Boyd, Jan Stienstra, and Michael Somos
for their useful comments and encouragement. Special thanks are due to James Wan
for helpful assistance in numerics and experimentation as well as for many
fruitful suggestions. We thank the anonymous referee for several valuable comments
which helped us to improve the exposition.
\end{acknowledgements}

\end{document}